\documentclass[seceqn,dvips]{arxbj}
\usepackage{euscript}
\usepackage{cursive,mathbh}  %uplgreek

\aid{0}
\volume{16}
\issue{4}
\pubyear{2010}
\firstpage{926}
\lastpage{952}
\doi{10.3150/09-BEJ244}

\makeatletter
\def\I{\mathbh{1}}
\def\l{\langle}
\def\dimh{\operatorname{dim}_{\mathrm{H}}}
\def\dimp{\operatorname{dim}_{\mathrm{P}}}
\def\Dim{\operatorname{Dim}}

\def\R{\mathbb{R}}
\def\C{\mathbb{C}}
\def\d{\mathbb{R}^d}

\def\a{\alpha}
\def\la{\lambda}
\def\ga{\gamma}
\def\ep{\varepsilon}
\def\eps{\varepsilon}

\def\Re{\operatorname{Re}}
\def\E{\mathbb{E}}
\def\P{\mathbb{P}}
\def\N{\mathbb{R}^N}

\newtheorem{proposition}{Proposition}[section]
\newtheorem{lemma}[proposition]{Lemma}
\newtheorem{theorem}[proposition]{Theorem}
\newtheorem{corollary}[proposition]{Corollary}

\newremark{remark}[proposition]{Remark}
\newremark{example}[proposition]{Example}

\makeatother

\begin{document}
\begin{frontmatter}

\title{Hausdorff and packing dimensions of the images of random fields}
\runtitle{Hausdorff and packing dimensions of the images}

\begin{aug}
\author[a]{\fnms{Narn-Rueih} \snm{Shieh}\thanksref{a}\ead[label=e1]{shiehnr@math.ntu.edu.tw}}
\and
\author[b]{\fnms{Yimin} \snm{Xiao}\thanksref{b}\ead[label=e2]{xiao@stt.msu.edu}\corref{}}
\runauthor{N.-R. Shieh and Y. Xiao}
\address[a]{Department of Mathematics, National Taiwan University,
Taipei, 10617, Taiwan.\\ \printead{e1}}
\address[b]{Department of Statistics and Probability, Michigan State
University, East Lansing, MI 48824, USA. \printead{e2}}
\end{aug}

% HISTORY:
\received{\smonth{4} \syear{2009}}
\revised{\smonth{9} \syear{2009}}

% ABSTRACT
%
\begin{abstract}
Let $X = \{X(t), t \in\mathbb{R}^N\}$ be a random field with values in
$\mathbb{R}^d$. For any finite Borel measure $\mu$ and analytic set $E
\subset\mathbb{R}^N$, the Hausdorff and packing dimensions of
the image measure $\mu_{X}$ and image set $X(E)$ are determined
under certain mild conditions.
These results are applicable to Gaussian random fields, self-similar
stable random fields with stationary increments, real harmonizable
fractional L\'evy fields and the Rosenblatt process.
\end{abstract}

% KEYWORDS
%

\pdfkeywords{Hausdorff dimension, images, packing dimension, packing dimension
profiles, real harmonizable fractional Levy motion, Rosenblatt process, self-similar stable random fields}

\end{frontmatter}

%s1 ###
\section{Introduction}

Fractal dimensions such as Hausdorff dimension, box-counting
dimension and packing dimension are very useful in characterizing
roughness or irregularity of stochastic processes and random fields
which, in turn, serve as stochastic models in various scientific
areas including image processing, hydrology, geostatistics and
spatial statistics. Many authors have studied the Hausdorff
dimension and exact Hausdorff measure of the image sets of Markov
processes and Gaussian random fields. We refer to Taylor (\citeyear{Taylor86a})
and Xiao (\citeyear{Xiao04}) for extensive surveys on results and techniques
for Markov processes, and to Adler (\citeyear{Adler81}) and Kahane (\citeyear{Kahane85}) for
results on Gaussian random fields.

Let $X = \{X(t), t \in\R^N\}$ be a random field with values
in $\R^d$, which will simply be called an $(N, d)$\emph{-random field}.
For any finite Borel measure $\mu$ on $\R^N$, the image measure
of $\mu$ under $X$ is defined by $\mu_{X}:=\mu\circ X^{-1}$.
Similarly, for every $E \subset\R^N$, the image set is
denoted by $X(E) = \{X(t), t \in E\} \subset\R^d$.
This paper is concerned with the Hausdorff and packing
dimensions of the image measures and image sets of random
fields which are, in a certain sense, comparable to a self-similar
process. Recall that $X=\{X(t), t \in\R^N\}$ is said to be
$H$\emph{-self-similar} if, for every constant $c > 0$, we have
%
%e1.1 ###
\begin{equation}\label{Eq:SS0}
\{X(c t), t \in\R^N \} \stackrel{d}{=}
\{c^H X(t), t \in\R^N \}
\end{equation}
and $X$ is said to have \emph{stationary increments} if, for every $h
\in
\R^N$,
%
%e1.2 ###
\begin{equation}\label{Eq:SI}
\{X(t + h) - X(h), t \in\R^N \} \stackrel{d}{=}
\{ X(t)- X(0), t \in\R^N \},
\end{equation}
where $\stackrel{d}{=}$ denotes equality of all
finite-dimensional distributions. If $X$ satisfies both
(\ref{Eq:SS0}) and (\ref{Eq:SI}), then it is called an $H$\emph{-SSSI
random field}. Samorodnitsky and Taqqu (\citeyear{ST94}) give a systematic
account of self-similar stable processes.
The main results of this paper show that the Hausdorff and
packing dimensions of the images of an $H$-SSSI random field
$X$ are determined by the self-similarity index $H$ and essentially
do not depend on the distributions of $X$.

An important example of an $H$-SSSI $(N, d)$-random field is fractional
Brownian motion $X = \{X (t), t \in\R^N\}$ of index $H$ $(0 < H < 1)$,
which is a centered Gaussian random field with the covariance function
$
\E[X_i(t)X_j(s)]=\frac{1}{2}\delta_{i,j} (\|s\|^{2H}+\|t\|^{2H}
-\|t-s\|^{2H} )$ for all $ s, t \in\R^N$,
where $\delta_{i, j} = 1$ if $i = j$ and $\delta_{i, j} =0$ otherwise.
It is well known (see Kahane (\citeyear{Kahane85}), Chapter 18) that for every
Borel set $E \subset\R^N$,
%
%e1.3 ###
\begin{equation}\label{Eq:BG60}
\dimh X(E) = \min \biggl\{d, \frac1 H \dimh E \biggr\}\qquad
\mathrm{a.s.},
\end{equation}
where $\dimh$ denotes the Hausdorff dimension.
On the other hand, Talagrand and Xiao (\citeyear{TX96}) proved that,
when $N > Hd$, the packing dimension analog of (\ref{Eq:BG60})
fails in general. Xiao (\citeyear{Xiao97a}) proved that
%
%e1.4 ###
\begin{equation} \label{Eq:Xiao97}
\dimp X(E) = \frac1 H \Dim_{Hd} E\qquad \mathrm{a.s.},
\end{equation}
where $\dimp$ denotes packing dimension and $\Dim_s E$ is
the packing dimension profile of $E$ defined
by Falconer and Howroyd (\citeyear{FH97}) (see Section \ref{Sec:Pre} for its definition).
Results (\ref{Eq:BG60}) and (\ref{Eq:Xiao97}) show that there
are significant differences between Hausdorff dimension and packing
dimension, and both dimensions are needed for characterizing the
fractal structures of $X(E)$.

There have been various efforts to extend (\ref{Eq:BG60}) to other
non-Markovian processes or random fields, but with only partial
success; see K\^ono (\citeyear{Kono86}), Lin and Xiao (\citeyear{LX94}),
Benassi, Cohen and Istas (\citeyear{BenassiCohen03}) and Xiao (\citeyear{Xiao07}). In order to
establish a Hausdorff
dimension result similar to (\ref{Eq:BG60}) for a random field $X$,
it is standard to determine upper and lower bounds for $\dimh X(E)$
separately. While the capacity argument (based on Frostman's
theorem) is useful for determining lower bounds, the methods
based on the classical covering argument for establishing an upper
bound for $\dimh X(E)$ are quite restrictive and usually require strong
conditions to be imposed on~$X$. As such, the aforementioned authors
have only considered random fields which either satisfy a uniform
H\"older condition of appropriate order on compact sets or
have at least the first moment. In particular, the existing
methods are not enough, even for determining $\dimh X([0,1]^N),$
when $X = \{X (t), t \in\R^N\}$ is a general stable random
field.

Given a random field $X= \{X (t), t \in\R^N\}$ and a Borel set $E
\subset\R^N$, it is usually more difficult to determine the packing
dimension of the image set $X(E)$. Recently, Khoshnevisan and Xiao
(\citeyear{KX08a}) and
Khoshnevisan, Schilling and Xiao (\citeyear{KSX09}) have solved the above
problem when $X = \{X(t), t \ge0\}$ is a L\'evy process in $\R^d$.
However, their method depends crucially on the strong Markov
property of L\'evy processes and cannot be applied directly
to random fields.

This paper is motivated by the need to develop methods for
determining the Hausdorff and packing dimensions of the image
measure $\mu_{X}$ and image set $X(E)$ under minimal conditions
on the random field $X$. By applying measure-theoretic
methods and the theory of packing dimension
profiles, we are able to solve the problems
for the Hausdorff and packing dimensions of the image
measure $\mu_{X}$ under mild conditions (namely, (C1)
and (C2) in Section \ref{Sec:measures}). The main results are Theorem
\ref{Th:Hdim1} and Theorem \ref{Th:main1}. When $X$ satisfies
certain uniform H\"older conditions, Theorems
\ref{Th:Hdim1} and \ref{Th:main1} can be applied directly to compute
the Hausdorff and packing dimensions of $X(E)$.
More generally, we also provide a method for determining the
Hausdorff dimension of $X(E)$ under conditions (C1) and (C2)
(see Theorem \ref{Th:dimXE2}). However, we have not been able
to solve the problem of determining $\dimp X(E)$ in general.

The rest of this paper is organized as follows. In Section \ref{Sec:Pre}, we
recall the definitions and some basic properties of
Hausdorff dimension, packing dimension and packing dimension
profiles of sets and Borel measures. In Section \ref{Sec:measures}, we determine
the Hausdorff and packing dimensions of the image measure
$\mu_{X}$ under general conditions (C1) and (C2). In Section \ref{Sec:set},
we study the Hausdorff and packing dimensions of the image set
$X(E)$, where $E \subset\R^N$ is an analytic set (i.e., $E$ is
a continuous image of the Baire space ${\mathbb N}^{\mathbb N}$ or,
equivalently, $E$ is a continuous image of a Borel set).
Section \ref{sec:examples} contains applications of the theorems in Sections \ref{Sec:measures} and \ref{Sec:set}
to SSSI stable random fields, real harmonizable fractional
L\'evy fields and the Rosenblatt process.

Throughout this paper, we will use $\langle x, y\rangle$ to
denote the inner product and $ \|\cdot\|$ to denote the
Euclidean norm in $\R^n$, no matter what the value of $n$ is.
For any $s, t \in\R^n$ such that $s_j < t_j$ ($j = 1,
\ldots, n$), $[s, t] = \prod^n_{j=1} [s_j, t_j]$ is called a
\emph{closed interval}.
We will use $K$ to denote an unspecified positive constant which
may differ from line to line. Specific constants in Section $i$
will be denoted by $K_{i,1}, K_{i,2},\ldots.$

%s2 ###
\section{Preliminaries}
\label{Sec:Pre}

In this section, we recall briefly the definitions and some basic
properties of Hausdorff dimension, packing dimension and packing
dimension profiles. More information on Hausdorff and packing
dimensions can be found in Falconer (\citeyear{Fal90}) and Mattila (\citeyear{Mattila95}).

%s2.1 ###
\subsection{Hausdorff dimension of sets and measures}

For any $\alpha> 0$, the $\alpha$-dimensional Hausdorff measure
of $E \subset\R^N$ is defined by
%
%e2.1 ###
\begin{equation}
\label{Eq:Hausdorff}
s^\alpha\mbox{-}m(E) = \lim_{\eps\to0}
\inf \Biggl\{ \sum_{i=1}^\infty(2 r_i)^\alpha\dvt E \subset
\bigcup_{i =1}^{\infty} B(x_i, r_i), r_i < \eps \Biggr\},
\end{equation}
where $B(x,r)= \{y \in\R^N\dvt |y-x|< r\}$. The Hausdorff
dimension of $E$ is defined as
$\dimh E = \inf \{ \a> 0\dvt \hbox{$s^{\a}$-$m$}(E) = 0 \}.
$
For a finite Borel measure $\mu$ on $\R^N$, its Hausdorff dimension
is defined by
$
\dimh\mu= \inf \{ \dimh E\dvt \mu(E) > 0 \mbox{ and }E \subset
\R^N\mbox{ is a Borel set} \}
$
and its upper Hausdorff dimension is defined by
$
\dimh^* \mu= \inf \{ \dimh E\dvt \mu({\R^N} \backslash E) = 0
\mbox{ and }E \subset\R^N\mbox{ is a Borel set} \}.
$
Hu and Taylor (\citeyear{HT94}) proved that
%
%e2.2 ###
\begin{eqnarray}\label{Eq:HuTaylor1}
\dimh\mu&=& \sup\biggl\{\beta> 0\dvt \limsup_{r \to0}
\frac{\mu(B(x, r) )} {r^{\beta}} = 0 \mbox{ for $\mu$-a.a. $
x \in\N$} \biggr\},
\\
\label{Eq:HuTaylor2}
\dimh^* \mu&=& \inf\biggl\{\beta> 0\dvt \limsup_{r \to0}
\frac{\mu(B(x, r) )} {r^{\beta}} > 0 \mbox{ for $\mu$-a.a. $
x \in\N$} \biggr\}.
\end{eqnarray}

The Hausdorff dimensions of an analytic set $E \subset\R^N$ and
finite Borel measures on $E$ are related by the following
identity (which can be verified by
(\ref{Eq:HuTaylor1}) and Frostman's lemma):
%
%e2.4 ###
\begin{equation}\label{Eq:setdim}
\dimh E = \sup \{ \dimh\mu\dvt \mu\in\EuScript{M}^+_c(E)
\},
\end{equation}
where $\EuScript{M}^+_c(E)$ denotes the family of all finite Borel
measures with compact support in $E$.

%s2.2 ###
\subsection{Packing dimension of sets and measures}

Packing dimension was introduced by Tricot (\citeyear{Tricot82}) as a dual
concept to Hausdorff dimension and has become a useful
tool for analyzing fractal sets and sample paths of
stochastic processes; see Taylor and Tricot (\citeyear{TT85}),
Taylor (\citeyear{Taylor86a}), Talagrand and Xiao (\citeyear{TX96}), Falconer and
Howroyd (\citeyear{FH97}), Howroyd (\citeyear{Howroyd01}), Xiao (\citeyear{Xiao97a}, \citeyear{Xiao04}, \citeyear{Xiao09}),
Khoshnevisan and Xiao (\citeyear{KX08a}, \citeyear{KX08b}), Khoshnevisan,
Schilling and Xiao (\citeyear{KSX09}) and the references therein
for more information.

For any $\alpha> 0$, the $\a$-dimensional packing measure
of $E \subset\R^N$ is defined as
\[
\hbox{$s^\alpha$-$p$} (E) = \inf \biggl\{\sum_n \hbox
{$\phi$-$P$}(E_n) \dvt E \subset\bigcup_n E_n \biggr\},
\]
where $s^\alpha$-$P$ is the set function on subsets of $\R^{N}$
defined by
\[
\hbox{$s^\alpha$-$P$}(E) = \lim_{\ep\to0} \sup \biggl\{ \sum_i
(2 r_i)^\alpha\dvt \overline{B}(x_i, r_i) \mbox{ are disjoint, } x_i
\in E, r_i < \ep \biggr\}.
\]
The packing dimension of $E$ is defined by
$
\dimp E = \inf\{ \a> 0\dvt \hbox{$s^{\a}$-$p$} (E) = 0
\}.
$
It is well known that $0 \le\dimh E \le\dimp E\le N$ for every
set $E\subset\R^N$.
%The packing dimension can also be defined through the upper
%box-counting dimension, see Tricot (1982) or Falconer (1990, p.45).

The packing dimension of a finite Borel measure $\mu$ on $\R^N$
is defined by $
\dimp\mu= \inf\{ \dimp E\dvt \mu(E) > 0 \mbox{ and }E \subset
\R^N\mbox{ is a Borel set}\}$
and the upper packing dimension of $\mu$ is defined by $
\dimp^* \mu= \inf\{ \dimp E\dvt \mu({\N} \backslash E) = 0
\mbox{ and }E \subset\R^N\mbox{ is a Borel set}\}.$
In analogy to (\ref{Eq:setdim}), Falconer and Howroyd (\citeyear{FH97})
proved, for every analytic set $E \subset\R^N$,
that
%
%e2.5 ###
\begin{equation}\label{Eq:setDim}
\dimp E = \sup \{ \dimp\mu\dvt \mu\in\EuScript{M}^+_c(E)
\}.
\end{equation}

%s2.3 ###
\subsection{Packing dimension profiles}

Next, we recall some aspects of the packing dimension profiles of
Falconer and Howroyd (\citeyear{FH97}) and Howroyd (\citeyear{Howroyd01}).
%They are not only useful for determining the
%packing and box-counting dimensions of projections, but also for
% studying random fractals; see Xiao (1997, 2009), Khoshnevisan and
%Xiao (2008b),
%Khoshnevisan, Schilling and Xiao (2009).
For a finite Borel measure $\mu$ on $\R^N$ and for any $s
> 0$, let
\[
%begin{equation}\label{Eq:F}
F_s^{\mu} (x, r) = \int_{\R^N} \psi_s \biggl(\frac{x-y} {r} \biggr)
\,\mathrm{d}\mu(y)
\]
%
%end{equation}
be the potential with respect to the kernel
$\psi_s(x)= \min\{1, \|x\|^{-s}\}, \forall x \in\R^N.$

Falconer and Howroyd (\citeyear{FH97})
defined the packing dimension profile and the upper packing
dimension profile of $\mu$ as
%
%e2.6 ###
\begin{equation}\label{Def:Dimprof}
\Dim_{s}\mu= \sup \biggl\{ \beta\ge0\dvt \liminf_{r \to0}
\frac{F_s^{\mu} (x, r)} {r^{\beta} } = 0 \mbox{ for
$\mu$-a.a. $ x \in\R^N$} \biggr\}
\end{equation}
and
%
%e2.7 ###
\begin{equation}\label{Def:Dimprof2}
\Dim^*_{s}\mu= \inf \biggl\{ \beta>0 \dvt \liminf_{r \to0}\frac{
F_s^{\mu} (x, r)} {r^{\beta} } > 0 \mbox{ for
$\mu$-a.a. $ x \in\R^N$} \biggr\},
\end{equation}
respectively. Further, they showed that
$
0 \le\Dim_{s} \mu\le\Dim^*_s \mu\le s
$
and, if $s \ge N$, then
%
%e2.8 ###
\begin{equation}\label{Eq:Rel2}
\Dim_{s} \mu= \dimp\mu,\qquad
\Dim^*_s \mu = \dimp^* \mu.
\end{equation}
%
%Note that the identities in (\ref{Eq:Rel2}) give equivalent
%characterizations of $\dimp\mu$ and $\dimp^* \mu$ in terms
%of the potential $F^{\mu}_N (x, r)$.

Motivated by (\ref{Eq:setDim}),
Falconer and Howroyd (\citeyear{FH97}) defined the $s$-dimensional packing
dimension profile of $E \subset\N$ by
%
%e2.9 ###
\begin{equation}\label{Eq:setDimprofile}
\Dim_{s} E = \sup \{ \Dim_{s} \mu\dvt \mu\in
\EuScript{M}^+_c(E) \}.
\end{equation}
It follows that
%
%e2.10 ###
\begin{equation} \label{Eq:setprofile2}
0 \le\Dim_{s} E \le s \quad\mbox{and}\quad \Dim_{s} E = \dimp E
\qquad\mbox{if } s \ge N.
\end{equation}
By the above definition, it can be verified (see Falconer and Howroyd (\citeyear{FH97}), page 286) that for every Borel
set $E\subset\N$ with $\dimh E = \dimp E$, we have
%
%e2.11 ###
\begin{equation}\label{Eq:setprofile3}
\Dim_{s} E = \min \{s, \dimp E \}.
\end{equation}
%
%Note that Eq. (\ref{Eq:setprofile3}) can also be proven by using
%a probabilistic argument. Namely, (\ref{Eq:setprofile3}) follows
%directly from (\ref{Eq:BG60})
%and (\ref{Eq:Xiao97}) by letting $s = Hd$.

The following lemma is a
consequence of Proposition 18 in Falconer and Howroyd (\citeyear{FH97}).

\begin{lemma}\label{Lem:Con}
Let $\mu$ be a finite Borel measure on $\R^N$ and $E \subset\R^N$
be bounded and non-empty. Let $\sigma\dvtx {\R}_+ \to[0, N]$ be any
one of the functions $\Dim_{s} \mu,$ $ \Dim^*_{s} \mu$ or $\Dim
_{s} E$ in $s$. Then $\sigma(s)$ is non-decreasing and continuous.
\end{lemma}

%s3 ###
\section{Hausdorff and packing dimensions of the image measures}
\label{Sec:measures}

Let $X = \{X(t), t \in\R^N\}$ be an $(N, d)$-random field defined
on some probability space $(\Omega, \mathcal{F}, \P)$. We assume
throughout this paper that $X$ is separable (i.e., there exists
a countable and dense set $T^* \subset\R^N$ and a zero probability
event $\Upsilon_0$ such that for every open set $F \subset\R^N$ and
closed set $G \subset\R^d$, the two events $\{\omega\dvt X(t, \omega)
\in G \mbox{ for all } t \in F\cap T^* \}$ and $\{\omega\dvt X(t, \omega)
\in G \mbox{ for all } t \in F \}$ differ from each other only by
a subset of $\Upsilon_0$; in this case, $T^*$ is called a \emph{separant}
for $X$) and $(t,\omega) \mapsto X(t, \omega)$
is $\mathcal{B}(\R^N)\times\mathcal{F}$-measurable, where $\mathcal{B}(\R^N)$
is the Borel $\sigma$-algebra of $\R^N$.

For any Borel measure $\mu$ on $\R^N $, the image measure $\mu_{X}$
of $\mu$ under $t \mapsto X(t)$ is
\[
\mu_{X}(B) := \mu\bigl\{ t \in\R^N \dvt X(t) \in B \bigr\}
\qquad\mbox{for all Borel sets } B \subset\d.
\]
In this section, we derive upper and lower bounds for the Hausdorff
and packing dimensions of the image measures of $X$, which rely,
respectively, on the following conditions (C1) and (C2). Analogous
problems for the image set $X(E)$ will be considered in Section \ref{Sec:set}.
\begin{enumerate}[(C2)]
\item[(C1)] There exist positive and finite constants $H_1$ and $\beta$
such that
%
%e3.1 ###
\begin{equation}\label{Eq:C1}
\P\Bigl\{\sup_{\|s- t\| \le h} \|X(s) - X(t) \| \ge h^{H_1}
u \Bigr\} \le K_{{3,1}} u^{- \beta}
\end{equation}
for all $ t \in\R^N$, $h \in(0, h_0)$ and $u \ge u_0$, where
$h_0$, $u_0$ and $K_{{3,1}}$ are positive constants.

\item[(C2)] There exists a positive constant $H_2$
such that for all $s, t \in\R^N$ and $r > 0$,
%
%e3.2 ###
\begin{equation}\label{Eq:C2}
\P\bigl\{ \|X(s) - X(t) \| \le\|s-t\|^{H_2} r\bigr\} \le
K_{{3,2}} \min \{1, r^d \} ,
\end{equation}
where $K_{{3, 2}} >0$ is a finite constant.
\end{enumerate}

\begin{remark}
Since (C1) and (C2) play essential roles in this
paper, we will now make some relevant remarks about them.
\begin{itemize}
\item Condition (C1) is a type of local maximal inequality and
is easier to verify when the random
field $X$ has a certain approximate self-similarity. For example,
if $X$ is $H_1$-self-similar, then condition (C1)
is satisfied whenever the tail probability of $\sup_{\|s- t\| \le
1} \|X(s) - X(t) \|$ decays no slower than a polynomial rate;
see Proposition \ref{Prop:tail} below and Section \ref{sec:examples}.
It can also be verified directly for Gaussian or more general
infinitely random fields by using large deviations techniques
without appealing to self-similarity.

\item There may be different pairs of $(H_1, \beta)$ for which (C1)
is satisfied. We note that the formulae for Hausdorff and
packing dimensions of the images do not depend on the constant
$\beta> 0$, it is $\sup\{H_1\dvt \mbox{(C1)  holds for some }
(H_1, \beta)\}$ that determines the best upper bounds for the
Hausdorff and packing dimensions of the image measures.

\item For every point $t\in\R^N$, the local H\"older exponent
of $X$ at $t$ is defined as
\[
\alpha_{X}(t) = \sup\biggl\{\gamma> 0\dvt \lim_{\|s-t\| \to0}
\frac{\|X(s)- X(t)\|} {\|s-t\|^\gamma} = 0 \biggr\}.
\]
Condition (C1) and the Borel--Cantelli lemma imply that
$\alpha_{X}(t) \ge H_1$ almost surely (see (\ref{Eq:LIL}) below).
However, (C1) does not even imply sample path continuity of $X$.

\item In Section \ref{Sec:set}, the following, slightly weaker, form of condition (C2)
will be sufficient:
\begin{enumerate}[(C2$'$)]
\item[(C2$'$)] There exist positive constants $H_2$ and $K_{{3, 2}}$
such that (\ref{Eq:C2}) holds for all $s, t \in\R^N$ satisfying
$\|t-s\| \le1$ and $r > 0$.
\end{enumerate}

\item Condition (C2) (or (C$2'$)) is satisfied if, for all $s, t
\in\R^N$ (or those satisfying $\|s-t\| \le1$), the random vector
$ (X(s)-
X(t) )/\|s-t\|^{H_2}$ has a density function which is uniformly
bounded in $s$ and $t$. As shown by Proposition \ref{Prop:C2} below,
(C2) is significantly weaker than the latter.
\end{itemize}
\end{remark}

The following proposition gives a simple sufficient condition for
an SSSI process $X = \{X(t), t \in\R\}$ to satisfy
condition (C1). More precise information can be obtained if
further distributional properties of $X$ are known; see Section \ref{sec:examples}.

\begin{proposition}\label{Prop:tail}
Let $X= \{X(t), t \in\R\}$ be a separable, $H$-SSSI process
with values in $\R^d$.
If there exist positive constants $\beta> 0$ and $K_{{3,3}}$
such that $H \beta> 1$ and
%
%e3.3 ###
\begin{equation}\label{Eq:X1tail}
\P\{ \|X(1)\| \ge u \} \le K_{{3,3}} u^{-\beta}\qquad
\forall u \ge1 ,
\end{equation}
then there exists a positive constant $K_{{3,4}}$ such that for
all $u \ge1$,
%
%e3.4 ###
\begin{equation} \label{Eq:SStail}
\P\Bigl\{\sup_{ t \in[0, 1]} \|X(t)\| \ge u \Bigr\} \le
K_{{3,4}} u^{-\beta}.
\end{equation}
In particular, condition \textup{(C1)} is satisfied with $H_1 = H$ and
the same $\beta$ as in (\ref{Eq:X1tail}).
\end{proposition}

\begin{pf}
Without loss of generality, we can assume $d=1$.
Since the self-similarity index $H > 0$, we have $X(0)= 0$ a.s.
Let $T^* = \{t_n, n \ge0\}$ be a separant for $X= \{X(t), t \in
[0, 1]\}$. We assume that $0= t_0 < t_1 < t_2 <\cdots< t_n < \cdots$.
For any $n \ge2$, consider the random variables $Y_k$ $(1 \le k
\le n)$ defined by $Y_k = X(t_k) - X(t_{k-1}).$
For $1 \le i < j \le n$, let $S_{i, j} = \sum_{k=i}^j Y_k$. By the
stationarity of increments and self-similarity of $X$ and
(\ref{Eq:X1tail}), we derive that for any $u \ge1$,
%
%e3.5 ###
\begin{equation}\label{Eq:MSS1}
\P\Biggl\{ \Biggl|\sum_{k=i}^j Y_k \Biggr|\ge u \Biggr\} =
\P\biggl\{ |X(1) |\ge\frac{u} {(t_j - t_{i-1})^H} \biggr\}
\le K_{{3,3}} u^{-\beta} (t_j - t_{i-1} )^{H \beta}.
\end{equation}
Thus, condition (3.4) of Theorem 3.2 of Moricz, Serfling and Stout
(\citeyear{MSS79}) is satisfied with $g(i, j) = t_j - t_{i-1}$, $\alpha= H
\beta$ and $\phi(t) = t^\beta$. It is easy to see that the
non-negative function $g(i, j)$ satisfies their condition (1.2)
(i.e., $g(i, j) \le g(i, j+1)$ and $g(i, j) + g(j+1, k) \le Q
g(i, k)$ for $1 \le i \le j < k \le n$) with $Q =1$. It therefore
follows from Theorem 3.2 of Moricz, Serfling and Stout (\citeyear{MSS79}) that
there exists a constant $K_{{3,4}}$ (independent of $n$) such that
for all $u \ge1$,
%
%e3.6 ###
\begin{equation}\label{Eq:MSS2}
\P\Bigl\{\max_{1 \le j \le n} |X(t_{j}) |\ge u \Bigr\} =
\P\Biggl\{\max_{1 \le j \le n} \Biggl|\sum_{k=1}^j Y_k \Biggr|\ge
u \Biggr\} \le K_{{3,4}} u^{-\beta}.
\end{equation}
Letting $n \to\infty$ yields (\ref{Eq:SStail}), which, in turn,
implies that (C1) holds for $H_1 = H$.
\end{pf}

%Theorem 3.2 of Moricz, Serfling and Stout (1982) requires $H
%holds under the condition $H \beta= 1$ is still open; see Moricz,
%Serfling and Stout (1982, p.1039). Even though
%Proposition \ref{Prop:tail} holds for any strictly stable L\'evy
%process $X$ with index $\a\in(0, 2]$ [in this case we have $H =
%1/\alpha$ and $\beta= \alpha$], it is not known whether Proposition
%If $H \beta< 1$, generally (\ref{Eq:X1tail}) does not imply
%(\ref{Eq:SStail}) as shown by the linear fractional stable processes
%with index $\a\in(0, 1)$.

Next, we provide a necessary and sufficient condition for an
$(N,d)$-random field $X=\{X(t), t \in\R^N\}$ to satisfy
condition (C2) (or (C$2'$)). For any $r > 0$, let
\[
\phi_r(x) = \prod_{j=1}^d \frac{1 - \cos(2r x_j)}{2 \curpi r x_j^2}
,\qquad x \in\R^d.
\]

\begin{proposition}\label{Prop:C2}
Let $X= \{X(t), t \in\R^N\}$ be a random field with values in $\R^d$.
Condition \textup{(C2)} (or (\textup{C}$2'$))
then holds if and only if there exists a positive constant $K_{{3,5}}$
such that
for all $r > 0$ and all $s, t \in\R^N$ $($or for those satisfying $\|
t-s\|\le1)$,
%
%e3.7 ###
\begin{equation}\label{Eq:C2-con}
\int_{\R^d} \phi_r(x) \E\bigl( \mathrm{e}^{\mathrm{i} \l x, X(t) - X(s)\rangle /\|t-s\|
^{H_2}} \bigr) \,\mathrm{d}x
\le K_{{3,5}} \min\{1, r^d\} .
\end{equation}
\end{proposition}

\begin{remark}
Since $\phi_r (x) = \mathrm{O}(\|x\|^{-2})$ as $\|x \| \to\infty$, condition
(\ref{Eq:C2-con})
is significantly weaker than assuming that $(X(t) - X(s))/\|t-s\|
^{H_2}$ has
a bounded density and can be applied conveniently to SSSI processes.
We mention that (\ref{Eq:C2-con}) is also weaker than the integrability
condition in Assumption 1 on page 269 of Benassi, Cohen and Istas
(\citeyear{BenassiCohen03}). It can be shown that Theorem 2.1 in Benassi, Cohen and Istas
(\citeyear{BenassiCohen03}) still holds under (\ref{Eq:C2-con}) and their Assumption 2.
\end{remark}

\begin{pf*}{Proof of Proposition \ref{Prop:C2}}
Note that for every
$r > 0$, the function $\phi_r(x)$
is non-negative and is in $L^1(\R^d)$. The Fourier transform of
$\phi_r$ is
\[
\widehat{\phi_r}(z) = \prod_{j=1}^d \biggl( 1 - \frac{|z_j|}{2r}
\biggr)^+\qquad \forall z\in\R^d .
\]
In the above, $a^+ :=\max(a ,0)$ for all $a\in\R$. Since
$z\in B(0 ,r)$ implies that $1-(2r)^{-1}|z_j| \ge\frac12$,
we have $\I_{B(0, r)}(z) \le2^d \widehat{\phi_r}(z)$
for all $z\in\R^d$. Here, and in the sequel, $\I_A$
denotes the indicator function (or random variable) of the set (or
event) $A$. By Fubini's theorem, we obtain
\begin{eqnarray*}
\P \{ \|X(s) - X(t) \| \le\|s-t\|^{H_2} r \} &\le& 2^d
\E\biggl[\widehat{\phi_r} \biggl( \frac{X(t) - X(s)} {\|t-s\|^{H_2}} \biggr) \biggr]\\
&=& 2^d \int_{\R^d} \phi_r(x) \E\bigl( \mathrm{e}^{\mathrm{i} \l x, X(t) - X(s)\rangle /\|
t-s\|^{H_2}} \bigr)\, \mathrm{d}x .
\end{eqnarray*}
Hence, (\ref{Eq:C2-con}) implies condition (C2). On the other hand, we have
$\widehat{\phi_r}(z) \le\I_{B(0, 2\sqrt{d} r)}(z)$ for all $z \in
\R^d$.
Consequently,
\[
\E\biggl[\widehat{\phi_r} \biggl( \frac{X(t) - X(s)} {\|t-s\|^{H_2}} \biggr) \biggr]
\le\P \bigl\{ \|X(s) - X(t) \| \le2\sqrt{d} \|s-t\|^{H_2} r \bigr\} .
\]
Therefore, condition (C2) implies
(\ref{Eq:C2-con}). This completes the proof.
\end{pf*}

%s3.1 ###
\subsection{Hausdorff dimensions of the image measures}
\label{Sec:measures-H}

First, we consider the upper bounds for the Hausdorff dimensions of the
image measure $\mu_{X}$.

\begin{proposition}\label{Prop:Hdimub1}
Let $X= \{X(t), t \in\N\}$ be a random field with values in $\R^d$.
If condition \textup{(C1)} is satisfied, then for every finite Borel measure
$\mu$ on $\N$,
%
%e3.8 ###
\begin{eqnarray}\label{Eq:HUpp1}
&&\dimh\mu_{X} \le\min \biggl\{d, \frac1 {H_1} \dimh\mu
\biggr\} \quad\mbox{and}\nonumber\\[-8pt]\\[-8pt]
&&\dimh^* \mu_{X} \le\min \biggl\{d, \frac1 {H_1} \dimh^* \mu
\biggr\}\qquad \mathrm{a.s.}\nonumber
\end{eqnarray}
\end{proposition}

\begin{pf}
Let $\lambda> 1/\beta$ be a constant.
For any fixed $s \in\R^N$ and
the sequence $h_n = 2^{-n}$ ($n \ge1$), it follows from condition
(C1) that for all integers $n \ge
\max\{\log(1/h_0), \frac1 {\log2} u_0^{1/\lambda} \}$,
\[
%begin{equation}\label{Eq:Tail2}
\P\Bigl\{\sup_{\|t -s\| \le2^{-n}} \|X(t) - X(s) \| \ge
2^{- H_1 n} (\log2^n)^{\la} \Bigr\} \le K n^{-\beta\la}.
\]
%
%end{equation}
Since $\sum_{n=1}^\infty n^{-\beta\la} < \infty$, the
Borel--Cantelli lemma implies that almost surely
%
%e3.9 ###
\begin{equation}\label{Eq:LIL}
\sup_{\|t -s\| \le2^{-n}} \|X(t) - X(s) \| \le(\log2)^\lambda 2^{- H_1
n} n^{\la} \qquad\forall n \ge n_0,
\end{equation}
where $n_0 = n_0(\omega, s)$ depends on $\omega$ and $s$. By
Fubini's theorem, we derive that, for any finite Borel measure
$\mu$ on $\N$, almost surely (\ref{Eq:LIL}) holds for $\mu$-a.a.
$s \in\R^N$.

We now fix an $\omega\in\Omega$ such that (\ref{Eq:LIL}) is valid
for $\mu$-a.a. $s \in\R^N$ and prove that both inequalities in
(\ref{Eq:HUpp1}) hold. In the sequel, $\omega$ will be suppressed.

To prove the first inequality in
(\ref{Eq:HUpp1}), since $\dimh\mu_{X} \le d$ holds
trivially, we only need to prove that
$\dimh\mu_{X} \le\frac1 {H_1} \dimh\mu$. Without loss of\vspace*{1pt}
generality, we assume $\dimh\mu_{X}>0$ and take any $\ga\in(0,
\dimh\mu_{X}).$ Then, by (\ref{Eq:HuTaylor1}), we have
%
%e3.10 ###
\begin{equation}\label{Eq:Hupp11}
\limsup_{r \to0} r^{- \ga} \int_{\d} \I_{\{ \|y - x\| \le r\}}
\,\mathrm{d}\mu_{X} (y) = 0 \qquad\mbox{for $\mu_{X}$-a.a. } x \in\R^d.
\end{equation}
Equivalently to (\ref{Eq:Hupp11}), we have
%
%e3.11 ###
\begin{equation}\label{Eq:211}
\limsup_{r \to0} r^{- \ga} \int_{\R^N} \I_{\{ \|X(t) - X(s)\|
\le r\}} \,\mathrm{d}\mu(t) = 0 \qquad\mbox{for $\mu$-a.a. } s \in\R^N.
\end{equation}

Let us fix $s \in\R^N$ such that both
(\ref{Eq:LIL}) and (\ref{Eq:211}) hold. For any $\varepsilon>0$,
we choose $n_1 \ge n_0$ such that $n^\la\le2^{\varepsilon n}$
for all $ n \ge n_1$. By (\ref{Eq:LIL}), we can write
\begin{eqnarray}\label{Eq:212}
\int_{\R^N} \I_{\{ \|X(t) - X(s)\|\le r\}} \,\mathrm{d}\mu(t)
&\ge&\sum_{n = n_1}^\infty\int_{2^{-n -1} \le\|t-s\|< 2^{-n}}
\I_{\{ \|X(t) - X(s)\| \le r \}}
\,\mathrm{d}\mu(t)\nonumber\\[-8pt]\\[-8pt]
&\ge&\int_{\|t-s\|< 2^{-n_1}} \I_{\{\|t-s\| \le
r^{1/(H_1-\ep)}\}} \,\mathrm{d}\mu(t) .\nonumber
\end{eqnarray}
Hence, we have
\begin{eqnarray}\label{Eq:213}
\int_{\R^N} \I_{\{\|t-s\| \le r^{1/(H_1 -\varepsilon) }\}} \,\mathrm{d}
\mu(t) &\le&
\int_{\R^N} \I_{\{ \|X(t) - X(s)\| \le r\}} \,\mathrm{d}\mu(t)
\nonumber\\[-8pt]\\[-8pt]
&&{} + \int_{\|t-s\|\ge 2^{-n_1}} \I_{\{\|t-s\| \le r^{1/(H_1
-\varepsilon) }\}} \,\mathrm{d}\mu(t) .\nonumber
\end{eqnarray}
For the last integral, we have
%
%e3.12 ###
\begin{equation}\label{Eq:214}
\lim_{r \to0} r^{-\ga} \int_{\|t-s\|\ge 2^{-n_1}}
\I_{\{\|t-s\| \le r^{1/(H_1 -\varepsilon) }\}} \,\mathrm{d}\mu(t) = 0
\end{equation}
because the indicator function takes the value 0 when $r>0$ is
sufficiently small.

It follows from (\ref{Eq:211}), (\ref{Eq:213}) and (\ref{Eq:214})
that with $r = \rho^{H_1-\ep}$,
\begin{eqnarray}\label{Eq:215}
\limsup_{\rho\to0} \rho^{-(H_1 - \ep) \ga} \int_{\R^N}
\I_{\{ \|t-s\| \le\rho\}} \,\mathrm{d}\mu(t)
& =& \limsup_{r \to0} r^{-\ga} \int_{\R^N} \I_{\{\|t-s\| \le
r^{1/(H_1 -\varepsilon)}\}} \,\mathrm{d}\mu(t)\nonumber\\[-8pt]\\[-8pt]
&\le&\limsup_{r \to0} r^{-\ga} \int_{\R^N}\I_{\{ \|X(t) -
X(s)\| \le r\}} \,\mathrm{d}\mu(t)= 0.\qquad\mbox{}\nonumber
\end{eqnarray}
We have thus proven that (\ref{Eq:215}) holds almost surely for
$\mu$-a.a. $s \in\R^N$. This implies that $\dimh\mu\ge(H_1 - \ep)
\ga$ almost surely. Since $\ep>0$ and $\ga< \dimh\mu_{X}$
are arbitrary, (\ref{Eq:HUpp1}) follows.

To prove the second inequality in (\ref{Eq:HUpp1}),
it is sufficient to show that
$\dimh^* \mu_{X} \le\frac1 {H_1} \dimh^* \mu$ a.s. Let $\omega
\in\Omega$ be fixed as above. We take an arbitrary $\beta>
\dimh^* \mu$. By (\ref{Eq:HuTaylor2}), we have
%
%e3.13 ###
\begin{equation} \label{Eq:226}
\limsup_{\rho\to0} \rho^{-\beta} \int_{\R^N}\I_{\{\|t-s\|\le
\rho\}} \,\mathrm{d} \mu(t)
> 0 \qquad\mbox{for $\mu$-a.a. $s \in\N$}.
\end{equation}

By using (\ref{Eq:212}), we derive that for $x = X(s)$,
\begin{eqnarray}\label{Eq:227}
\int_{\R^d}\I_{\{\|y-x\| \le r\}} \,\mathrm{d} \mu_{X}(y) &\ge&
\int_{\R^N} \I_{\{ \|t-s\| \le r^{1/(H_1 - \varepsilon)}\}} \,\mathrm{d}\mu
(t)\nonumber\\[-8pt]\\[-8pt]
&&{} - \int_{\|t-s\|\ge 2^{-n_1}} \I_{\{ \|t-s\| \le r^{1/(H_1-
\varepsilon)}\}} \,\mathrm{d}\mu(t) .\nonumber
\end{eqnarray}
It follows from (\ref{Eq:227}), (\ref{Eq:214}) and (\ref{Eq:226})
that
%
%e3.14 ###
\begin{equation}\label{Eq:229}
\limsup_{r \to0} \frac{
\int_{\R^d}\I_{\{\|y-x\| \le r\}} \,\mathrm{d} \mu_{X}(y)}{r^{ \beta/(H_1 -
\ep)}}
\ge\limsup_{r \to0}\frac{\int_{\R^N}
\I_{\{ \|t-s\| \le r^{1/(H_1 - \varepsilon)}\}} \,\mathrm{d}\mu(t)}
{r^{ \beta/(H_1 - \ep)}} > 0
\end{equation}
for all $s \in\R^N$ that satisfy (\ref{Eq:226}). This implies that $
\dimh^*\mu_{X} \le{\beta}/{(H_1-\ep)}$ a.s. Letting $\ep
\downarrow0$ and $\beta\downarrow\dimh^* \mu$ yields the
second inequality in (\ref{Eq:HUpp1}). This completes the proof
of Proposition~\ref{Prop:Hdimub1}.\mbox{\quad}
\end{pf}

\begin{remark} \label{Re:Hdimub1}
Note that in (\ref{Eq:HUpp1}), the
exceptional null probability events depend on $\mu$. For several purposes,
it is more useful to have a single exceptional null probability event
$\Omega_0$ such that, for all $\omega\notin\Omega_0$, both
inequalities in
(\ref{Eq:HUpp1}) hold \emph{simultaneously} for
all finite Borel measures $\mu$ on $\R^N$. By slightly modifying the
proof of Proposition \ref{Prop:Hdimub1} (see (\ref{Eq:212})), one
can show that this is indeed true if, for every $\ep> 0$ and every
compact interval $I$, the sample function $X(t)$ satisfies almost surely
a uniform H\"older condition of order $H_1-\ep$ on $I$.
\end{remark}

Next, we show that condition (C2) determines lower bounds for the
Hausdorff dimensions of the image measures of the random field
$X$.

\begin{proposition}\label{Prop:Hdimlb1}
Let $X=\{X(t), t \in\R^N \}$ be an $(N, d)$-random field
satisfying condition \textup{(C2)}. Then, for every finite Borel measure
$\mu$ on $ \R^N$,
%
%e3.15 ###
\begin{equation}\label{Eq:HLow1}
\dimh\mu_{X} \ge\min \biggl\{d, \frac1 {H_2} \dimh\mu
\biggr\} \quad\mbox{and}\quad
\dimh^* \mu_{X} \ge\min \biggl\{d, \frac1 {H_2} \dimh^* \mu
\biggr\}\qquad \mathrm{a.s.}
\end{equation}
\end{proposition}

\begin{pf}
In order to prove the first inequality in (\ref{Eq:HLow1}),
we fix any constants $0 <
\gamma< \gamma' <\min\{d, \frac1 {H_2} \dimh\mu \}$.
Since $\dimh\mu> \ga' H_2$, it follows from (\ref{Eq:HuTaylor1})
that
%
%e3.16 ###
\begin{equation}\label{Eq:HLow3}
\limsup_{r \to0} \frac{\mu(B(s, r) )} {r^{\ga' H_2}} = 0
\qquad\mbox{for $\mu$-a.a. $s \in\R^N$}.
\end{equation}

Let $s \in\R^N$ be a fixed point such that (\ref{Eq:HLow3}) holds.
By (C2), we derive
\begin{eqnarray}\label{Eq:HLow4}
\E \mu_{X} (B(X(s), r) ) &=& \int_{\R^N} \P\bigl(\|X(t) -
X(s)\| \le r \bigr) \mu(\mathrm{d}t)\nonumber\\[-8pt]\\[-8pt]
&\le& K_{{3,2}} \mu(B(s, r^{1/H_2}) ) + K_{{3,2}}
\int_{\|t -s \| > r^{1/H_2}} \biggl(\frac r {\|t-s\|^{H_2}} \biggr)^d\mu(\mathrm{d}t) .\nonumber
\end{eqnarray}
Let $\kappa$ be the image measure of $\mu$ under the mapping $t
\mapsto\|t-s\|$ from $\R^N$ to $\R_+$. Then, by using
an integration-by-parts formula and (\ref{Eq:HLow3}), we have
%
%e3.17 ###
\begin{eqnarray}\label{Eq:HLow5}
\int_{\|t -s \| > r^{1/H_2}} \biggl(\frac r {\|t-s\|^{H_2}} \biggr)^d
\mu(\mathrm{d}t) &=& \int_{r^{1/H_2}}^\infty\frac{r^d} {\rho^{H_2 d}}
\kappa(\mathrm{d}\rho)\nonumber\\
&\le& H_2 d \int_{r^{1/H_2}}^\infty\frac{r^d} {\rho^{H_2 d +1}}
\mu(B(s, \rho) ) \,\mathrm{d}\rho\\
&\le& K r^{\ga'}\nonumber
\end{eqnarray}
for all $r>0$ small enough, where the last inequality follows from
(\ref{Eq:HLow3}) and the fact that $\gamma' < d$. Combining (\ref{Eq:HLow4})
and (\ref{Eq:HLow5}), we see that
$\E \mu_{X} (B(X(s), r) ) \le K r^{\ga'}$ for $r>0$ small.
This, and the Markov inequality, imply that for all $n$ large enough,
\[
%begin{equation}\label{Eq:HLow7}
\P\bigl(\mu_{X} (B(X(s), 2^{-n}) ) \ge2^{-n \ga} \bigr) \le
K 2^{-n(\ga'-\ga)}.
\]
%
%end{equation}
It follows from the Borel--Cantelli lemma that a.s. $\mu_{X} (B(X(s),
2^{-n}) ) < 2^{-n \ga}$ for all $n$ large
enough. It should be clear the above implies that for all $0 <
\gamma< \min\{d, \frac1 {H_2} \dimh\mu \}$,
\[
\limsup_{r \to0} \frac{\mu_{X} (B(x, r) )} {r^{\ga}} = 0
\qquad\mbox{for $\mu_{X}$-a.a. $x \in\R^d$}
\]
almost surely. Thus, $\dimh\mu_{X} \ge\ga$ a.s., and (\ref{Eq:HLow1})
follows from the arbitrariness of $\ga$.

To prove the second inequality in (\ref{Eq:HLow1}), let $0 < \gamma<
\ga'< \min \{d, \frac1 {H_2} \dimh^* \mu \}$. By
(\ref{Eq:HuTaylor2}), there exists a Borel set $A \subset
\R^N$ such that $\mu(A) > 0$ and
$ \limsup_{r \to0} r^{-\ga' H_2} \mu(B(s, r) ) =
0$ for all $ s \in A.$ The proof above shows that a.s. $\limsup_{r \to
0} r^{-\ga} \mu_{X} (B(x, r) ) = 0$
for all $x \in X(A)$.
Since $\mu_{X} (X(A) )>0$ a.s., we derive
$\dimh^* \mu_{X} \ge\ga$ a.s. and the proof is completed.
\end{pf}

Combining Propositions \ref{Prop:Hdimub1} and \ref{Prop:Hdimlb1}, we
have the following theorem, whose proof is omitted.
\begin{theorem}\label{Th:Hdim1}
Let $X=\{X(t), t \in\R^N \}$ be an $(N, d)$-random field and let
$H$ be a positive constant. If, for every $\ep> 0$, $X$ satisfies
condition \textup{(C1)} with $H_1 = H-\ep$, some $\beta= \beta(\ep)>0$
and \textup{(C2)} with $H_2=H +\ep$,
then for every finite Borel measure $\mu$ on $\R^N$,
%
%e3.18 ###
\begin{equation}\label{Eq:H-dim1}
\dimh\mu_{X} = \min \biggl\{d, \frac1 {H} \dimh\mu \biggr\}\quad
\mbox{and}\quad
\dimh^* \mu_{X} = \min \biggl\{d, \frac1 {H} \dimh^* \mu
\biggr\}\qquad \mathrm{a.s.}
\end{equation}
\end{theorem}

%s3.2 ###
\subsection{Packing dimensions of the image measures}
\label{Sec:measures-P}

We now study the problem of determining the packing dimensions $\dimp
\mu_{X}$
and $\dimp^* \mu_{X}$.
The following upper bounds for the image measures are proved by Schilling
and Xiao (\citeyear{SchXiao05}).

\begin{proposition}\label{Lem:UppBD}
Let $X= \{X(t), t \in\N\}$ be a random field with values in $\R^d$.
If condition \textup{(C1)} is satisfied, then for every finite Borel measure
$\mu$ on $\N$,
%
%e3.19 ###
\begin{equation}\label{Eq:Upp1}
\dimp\mu_{X} \le\frac1 {H_1} \Dim_{{H_1 d}} \mu
\quad\mbox{and}\quad
\dimp^* \mu_{X} \le\frac1 {H_1} \Dim^*_{{H_1 d}} \mu
\qquad\mathrm{a.s.}
\end{equation}
\end{proposition}

Similarly to Remark \ref{Re:Hdimub1}, we have the following.
\begin{remark}\label{Re:UppBD}
If, for every $\ep> 0$ and every compact interval $I \subset\R^N$,
$X(t)$ satisfies almost surely a uniform H\"older condition
of order $H_1 - \ep$ on $I$, then almost surely both inequalities in
(\ref{Eq:Upp1}) hold for all finite Borel measures $\mu$ on $\N$.
\end{remark}

For the lower bounds of packing dimensions, we have
the following proposition.
\begin{proposition}\label{Prop:LowerBD}
Let $X=\{X(t), t \in\R^N \}$ be an $(N, d)$-random field
satisfying condition \textup{(C2)}. Then, for every finite Borel measure
$\mu$ on $ \R^N$,
%
%e3.20 ###
\begin{equation} \label{Eq:mDim1}
\dimp \mu_X \ge\frac1 {H_2} \Dim_{{H_2 d}} \mu
\quad\mbox{and}\quad
\dimp^* \mu_X \ge\frac1
{H_2} \Dim^* _{{H_2 d}} \mu \qquad\mathrm{a.s.}
\end{equation}
\end{proposition}

\begin{pf}
We only prove the first inequality in
(\ref{Eq:mDim1}); the proof of the second one is similar.
We may, and will, assume that $\Dim_{{H_2 d}} \mu> 0$.
For fixed $s \in\N$, Fubini's
theorem implies that
%
%e3.21 ###
\begin{equation}\label{Eq:3lower}
{\E} F^{\mu_{X}}_d (X(s), r ) = \int_{\N} {\E} \min\{1, r^d \|X(t) - X(s)\|^{-d} \} \,\mathrm{d}\mu(t).
\end{equation}
The integrand in (\ref{Eq:3lower}) can be written as
\begin{eqnarray}\label{Eq:3lower2}
&&\E\min \{1, r^d \|X(t) - X(s)\|^{-d} \}\nonumber\\[-8pt]\\[-8pt]
&&\quad = {\P} \{ \|X(t) - X(s)\| \le r \} + {\E} \bigl\{r^d
\|X(t) - X(s)\|^{-d} \cdot\I_{\{ \|X(t) - X(s)\| \ge r\}} \bigr\} .\nonumber
\end{eqnarray}
By condition (C2), we obtain that for all $ s, t \in\R^N$ and $r
> 0$,
%
%e3.22 ###
\begin{equation} \label{Eq:3lower2a}
{\P} \{ \|X(t) - X(s)\| \le r \} \le K_{{3, 2}}
\min\biggl\{1, \frac{r^d} {\|t - s\|^{H_2 d}} \biggr\}.
\end{equation}
Denote the distribution of $X(t) - X(s)$ by $\Gamma_{s,
t}(\cdot)$. Let $\nu$ be the image measure of $\Gamma_{s,
t}(\cdot)$ under the mapping $T\dvtx z \mapsto\|z\|$ from $\R^d$ to
$\R_+$. The last term in (\ref{Eq:3lower2}) can then be written
as
\begin{eqnarray}\label{Eq:3lower3}
\int_{\R^d} \frac{r^d} {\|z\|^d} \I_{\{\|z\|\ge r\}}
\Gamma_{s,t}(\mathrm{d}z) &=& \int_{r}^\infty\frac{r^d}
{\rho^d} \nu(\mathrm{d}\rho)\nonumber\\[-8pt]\\[-8pt]
&\le& d \int_{r}^\infty\frac{ r^d}
{\rho^{d+1}} \P\{\|X(t) - X(s)\| \le\rho\} \,\mathrm{d}\rho,\nonumber
\end{eqnarray}
where the last inequality follows from an integration-by-parts
formula.

By (\ref{Eq:3lower2a}) and (\ref{Eq:3lower3}), we derive
that the last term in (\ref{Eq:3lower2}) can be bounded by
a constant multiple of
%
%e3.23 ###
\begin{eqnarray}\label{Eq:3lower4}
&&\int_{r}^\infty\frac{r^d}
{\rho^{d+1}} \min \biggl\{1, \frac{\rho^d}
{\|t-s\|^{H_2 d}} \biggr\} \,\mathrm{d}\rho\nonumber\\[-8pt]\\[-8pt]
&&\quad\le \cases{
K , &\quad if $ r \ge\|t-s\|^{H_2} $,\vspace*{2pt}\cr
\displaystyle\frac{K r^d} {\|t-s\|^{H_2d}} \log \biggl(\frac
{\|t-s\|^{H_2}}
{r} \biggr) , &\quad if $ r < \|t-s\|^{H_2} $.
}\nonumber
\end{eqnarray}
It follows from (\ref{Eq:3lower2}), (\ref{Eq:3lower2a}),
(\ref{Eq:3lower3}) and (\ref{Eq:3lower4}) that for any $0 < \ep<
1$ and $s, t \in\R^N$,
%
%e3.24 ###
\begin{equation}\label{Eq:3lower5}
{\E} \min\{1, r^d \|X(t) - X(s)\|^{-d}\} \le K_{{3, 6}} \min
\biggl\{1, \frac{r^{d-\ep}} {\|t- s\|^{H_{2}(d -\ep)}} \biggr\}.
\end{equation}

For any $\ga\in(0, \Dim_{{H_2 d}}\mu),$ by Lemma
\ref{Lem:Con}, there exists $\ep> 0$ such that $\ga< \Dim_{H_2
(d- \ep)}\mu.$ It follows from (\ref{Def:Dimprof}) that
%
%e3.25 ###
\begin{equation} \label{Eq:3lower6}
\liminf_{r \to0} r^{- {\ga}/{H_2}} \int_{\R_+}
\min\biggl\{1, \frac{r^{d-\ep}} {\|t- s\|^{H_2(d -\ep)}} \biggr\}
\,\mathrm{d}\mu(t) = 0 \qquad\mbox{for $\mu$-a.a. $s \in\R^N$}.
\end{equation}

By (\ref{Eq:3lower}), (\ref{Eq:3lower5}), (\ref{Eq:3lower6}) and
Fatou's lemma, we have that for $\mu$-a.a. $s\in\R^N$,
\begin{eqnarray}\label{Eq:3lower7}
&&{\E} \Bigl(\liminf_{r \to0} r^{- {\ga}/{H_2}}
F^{\mu_{X}}_d (X(s), r ) \Bigr)\nonumber\\[-8pt]\\[-8pt]
&&\quad\le K_{{3, 6}} \liminf_{r \to0} r^{- {\ga}/{H_2}}
\int_{\R^N} \min\biggl\{1, \frac{r^{d-\ep}} {\|t- s\|^{H_2(d
-\ep)}} \biggr\} \,\mathrm{d}\mu(t)= 0 .\nonumber
\end{eqnarray}
By using Fubini's theorem again, we see that almost surely
\[
\liminf_{r \to0} r^{-{\ga}/H_2} F^{\mu_{X}}_d (X(s),
r ) = 0 \qquad\mbox{for $\mu$-a.a. $s \in\R^N$}.
\]
Hence, $ \dimp\mu_{X} \ge\frac{\ga} {H_2}$ a.s. Since $\ga$ can
be arbitrarily close to $\Dim_{{H_2 d}} \mu$, we obtain (\ref{Eq:mDim1}).
\end{pf}

The following is a direct consequence of Propositions
\ref{Lem:UppBD} and \ref{Prop:LowerBD}.

\begin{theorem}\label{Th:main1}
Let $X=\{X(t), t \in\R^N \}$ be an $(N, d)$-random field and let
$H$ be a positive constant. If, for every $\ep> 0$, $X$ satisfies
condition \textup{(C1)} with $H_1 = H-\ep$, some $\beta= \beta(\ep)>0$ and
\textup{(C2)} with $H_2=H +\ep$, then for every finite Borel measure
$\mu$ on $ \R^N$,
%
%e3.26 ###
\begin{equation} \label{Eq:mDim3}
\dimp \mu_X = \frac1 {H} \Dim_{{H d}} \mu \quad\mbox{and}\quad
\dimp^* \mu_X = \frac1 {H} \Dim^*_{{H d}} \mu
\qquad\mathrm{a.s.}
\end{equation}
\end{theorem}

%s4 ###
\section{Hausdorff and packing dimensions of the image sets}
\label{Sec:set}

We now consider the Hausdorff and packing dimensions of
the image set $X(E)$. We will see that general lower bounds for $\dimh X(E)$
and $\dimp X(E)$ can be derived from the results in Section
\ref{Sec:measures} by using a measure theoretic method. For random
fields which satisfy uniform H\"older conditions on compact intervals, the
upper bounds for $\dimh X(E)$ and $\dimp X(E)$ can also be easily obtained.
However, it is difficult to obtain upper bounds for
$\dimh X(E)$ and $\dimp X(E)$ under condition (C1) alone. We
have only been able to provide a partial result on determining
the upper bound for $\dimh X(E)$. The analogous problem for $\dimp X(E)$
remains open.

We will need the following lemmas. Lemma \ref{Lem:Edgar} is
from Lubin (\citeyear{Lubin74}), which is more general than Theorem 1.20
in Mattila (\citeyear{Mattila95}).

\begin{lemma}\label{Lem:Edgar}
Let $E \subset\R^N$ be an analytic set and let $f\dvtx \R^N \to\R^d$
be a Borel function. If $\nu$ is a finite Borel measure on $\R^d$
with support in $f(E),$ then $\nu= \mu_{f}$ for some $\mu\in
{\EuScript
M}^+_c(E)$.
\end{lemma}

\begin{lemma}\label{Lem:setdim}
Let $E \subset\R^N$ be an analytic set. Then, for all Borel
measurable functions $f: \R^N \to\R^d$, we have
%
%e4.2 ###
%e4.1 ###
\begin{eqnarray} \label{Eq:41}
\dimh f(E) &=& \sup\{ \dimh\mu_{f} \dvtx \mu\in{\EuScript
M}^+_c(E) \} ,
\\ \label{Eq:42}
\dimp f(E) &=& \sup\{ \dimp\mu_{f} \dvtx \mu\in{\EuScript
M}^+_c(E) \} .
\end{eqnarray}
\end{lemma}

\begin{pf}
Denote the right-hand side of (\ref{Eq:41}) by
$\gamma_{E}$. By (\ref{Eq:setdim}), we get $\dimh f(E)
\ge\ga_{E}$. Next, for any $\nu\in{\EuScript
M}_c^+(f(E))$, Lemma \ref{Lem:Edgar} implies that $\nu= \mu_{f}$
for some $\mu\in{\EuScript M}^+_c(E)$. This and
(\ref{Eq:setdim}) together imply $\dimh f(E) \le\ga_{E}$. Hence,
(\ref{Eq:41}) is proved.
The proof of (\ref{Eq:42}) is similar and is therefore omitted.
\end{pf}

We first consider the lower bounds for the Hausdorff and packing
dimensions of $X(E)$.
\begin{proposition}\label{Prop:LowerBD2}
Let $X=\{X(t), t \in\R^N \}$ be an $(N, d)$-random field that
satisfies condition (\textup{C}$2'$). Then, for every analytic set $E\subset
\R^N$,
%
%e4.3 ###
\begin{equation}\label{Eq:setLB1}
\dimh X(E) \ge\min\biggl\{d, \frac1 {H_2} \dimh E \biggr\}
\quad\mbox{and}\quad
\dimp X(E) \ge\frac1 {H_2} \Dim_{{H_2 d}} E \qquad\mathrm{a.s.}
\end{equation}
\end{proposition}
\begin{pf}
Since both $\dimh$ and $\dimp$ are $\sigma$-stable
(see Falconer (\citeyear{Fal90})), we may, and will,
assume that the diameter of $E$ is at most 1. Hence, condition
(C2$'$) will be enough to prove (\ref{Eq:setLB1}).

Let us prove the first inequality in (\ref{Eq:setLB1}). It follows from
(\ref{Eq:setdim}) that for any $0 < \ga< \dimh E$, there exists
a $\mu\in\EuScript{M}^+_c(E)$ such that
$\dimh\mu\ge\ga$. By Proposition \ref{Prop:Hdimlb1} (which holds for
any finite Borel measure whose support has diameter $\le1)$, we have\vspace*{-1.5pt}
$\dimh\mu_{X} \ge\min\{d, \frac1 {H_2} \dimh\mu
\}$ a.s. This and (\ref{Eq:41}) together imply that $\dimh X(E) \ge
\min\{d, \frac1 {H_2} {\ga}\}$ a.s. Since $\ga< \dimh
E$ is arbitrary, the desired inequality follows.

Next, we prove the second inequality in (\ref{Eq:setLB1}).
Note that for any $0 < \gamma<
\frac1 {H_2} \Dim_{{H_2 d}} E$, by (\ref{Eq:setDimprofile}),
there exists a Borel measure $\mu\in{\EuScript M}^+_c(E)$ such
that $H_2 \ga< \Dim_{{H_2 d}} \mu.$ It follows from
(\ref{Eq:mDim1}) that $ \dimp\mu_{X} > \ga$ a.s. Hence, by Lemma
\ref{Lem:setdim}, we have $ \dimp X(E) > \ga$ a.s., which, in turn,
implies that $\dimp X(E) \ge\frac1 {H_2} \Dim_{{H_2 d}} E$ a.s.
The proof is therefore completed.
\end{pf}

The following proposition gives upper bounds for
$\dimh X(E)$ and $\dimp X(E)$.

\begin{proposition}\label{Prop:setUpp1}
Let $X= \{X(t), t \in\N\}$ be an ($N, d$)-random field. If
for every $\ep> 0$, $X$ satisfies a uniform H\"older
condition of order $H_1-\ep$ on all compact
intervals of $\R^N$ almost surely, then, for all
analytic sets $E \subset\R^N$,
%
%e4.4 ###
\begin{equation}\label{Eq:setUpp1}
\dimh X(E) \le\min\biggl\{d, \frac1 {H_1} \dim E \biggr\}
\quad\mbox{and}\quad
\dimp X(E) \le\frac1 {H_1} \Dim_{{H_1 d}} E\qquad \mathrm{a.s.}
\end{equation}
\end{proposition}

\begin{pf}
Both inequalities in (\ref{Eq:setUpp1})
follow from Remarks \ref{Re:Hdimub1}, \ref{Re:UppBD} and Lemma
\ref{Lem:setdim}.
\end{pf}

Combining Propositions \ref{Prop:LowerBD2} and
\ref{Prop:setUpp1} yields the following theorem.

\begin{theorem}\label{Th:main2a}
Let $X=\{X(t), t \in\R^N \}$ be an $(N, d)$-random field and
let $H \in(0, 1]$ be a constant. If, for every $\ep> 0$, $X$ satisfies
a uniform H\"older condition of order $H-\ep$ on all compact
intervals of $\R^N$ and condition (\textup{C}$2'$) with $H_2=H +\ep$, then,
for all analytic sets $E\subset\R^N$,
%
%e4.5 ###
\begin{equation} \label{Eq:setResult4a}
\dimh X(E) = \min\biggl\{d, \frac1 {H} \dimh E \biggr\}\quad
\mbox{and}\quad
\dimp X(E) = \frac1 {H} \Dim_{Hd} E \qquad\mathrm{a.s.}
\end{equation}
\end{theorem}

It is often desirable to compute $\dimp X(E)$ in terms of $\dimp
E$. The following is the packing dimension analog of
(\ref{Eq:BG60}). Note that if $N > Hd$, then
the conclusion of Corollary \ref{Co:46} does not hold in general;
see Talagrand and Xiao (\citeyear{TX96}). In this sense, it is the best
possible result of this kind.

\begin{corollary}\label{Co:46}
Let $X=\{X(t), t \in\R^N \}$ and $E \subset\R^N$ be as in
Theorem \ref{Th:main2a}. If either $N \le H d$ or $E$
satisfies $\dimh E = \dimp E,$ then
$\dimp X(E) = \min\{d, \frac{1} {H} \dimp E \}$
a.s.
\end{corollary}

\begin{pf}
If $N \le Hd$, then (\ref{Eq:setprofile2}) implies that for every
analytic set $E \subset\N$, $\Dim_{{H d}} E = \dimp E.$ Hence, Theorem
\ref{Th:main2a} yields $\dimp X(E) = \frac1 {H} \dimp E$ a.s.,
as desired. On the other hand, if an analytic set $E \subset\N$
satisfies $\dimh E = \dimp E,$ then (\ref{Eq:setprofile3})
implies that $\Dim_{{H d}}E = \min\{ H d, \dimp E \}$.
Hence, again, the conclusion follows from Theorem \ref{Th:main2a}.
\end{pf}

Since many random fields do not have continuous sample functions and, even
if they do, it is known that $\dimh X(E)$ is not
determined by the exponent of uniform
modulus of continuity (a~typical example being linear fractional stable
motion -- see Example \ref{Ex:lfs} below), there have been various efforts
to remove the uniform H\"older condition. However, except for Markov
processes or random fields with certain Markov structure, no satisfactory
method has been developed. The main difficulty lies in deriving a sharp
upper bound for $\dimh X(E)$.

In the following, we derive an upper bound for $\dimh X(E)$ under
condition (C1). This method is partially motivated by an argument
in Schilling (\citeyear{Sch98}) for Feller processes generated by
pseudo-differential operators and, as far as we know, is more
general than the existing methods in the literature.

\begin{lemma}\label{Lem:Momentest}
Let $X=\{X(t), t \in\R^N\}$ be a random field with values in $\R^d$.
If condition \textup{(C1)} holds for $H_1>0$ and $\beta> 0$, then, for
all $t \in\R^N$, $h > 0$ and $\gamma> 0$,
%
%e4.6 ###
\begin{equation}\label{Eq:Momentest}
\E\bigl(D(t, h)^\gamma \mathrm{e}^{- D(t, h)} \bigr) \le K_{{4,1}}
h^{H_1(\gamma\wedge\beta)} ,
\end{equation}
where $D(t, h) = \sup_{\|s-t\|\le h} \|X(s) - X(t)\|$ and $K_{{4,1}}$ is
a constant independent of $t$ and $h$.
\end{lemma}

\begin{pf}
We write
\begin{eqnarray}\label{Eq:Momentest2}
\E\bigl(D(t, h)^\gamma \mathrm{e}^{- D(t, h)} \bigr) &=&
\int_0^\infty u^{\gamma- 1} \mathrm{e}^{- u} (\gamma- u)
\P\{D(t, h) > u \} \,\mathrm{d}u\nonumber\\[-8pt]\\[-8pt]
&\le& K \int_0^\gamma u^{\gamma- 1} \mathrm{e}^{- u} (\gamma- u) \min\{1,
(h^{-H_1} u)^{-\beta}\} \,\mathrm{d}u ,\nonumber
\end{eqnarray}
where the inequality follows from (C1). It is elementary
to verify that, up to a constant, the last integral is
bounded by
%
%e4.7 ###
\begin{equation}\label{Eq:Momentest3}
\int_0^{h^{H_1}} u^{\gamma- 1} \,\mathrm{d}u + h^{H_1 \beta}
\int_{h^{H_1}}^\gamma u^{\gamma- \beta- 1} (\gamma - u)
\,\mathrm{d}u \le K_{{4,1}} h^{H_1(\gamma\wedge\beta)} .
\end{equation}
This proves (\ref{Eq:Momentest}).
\end{pf}

\begin{proposition}\label{Prop:setdimUppBD}
Let $X=\{X(t), t \in\R^N\}$ be a random
field with values in $\R^d$. Suppose that the sample function
of $X$ is a.s. bounded on all compact subsets of $\R^N$. If
condition \textup{(C1)} holds for $H_1>0$ and $\beta> 0$, then, for every
analytic set $E \subset\R^N$ that satisfies $\dimh E < \beta H_1$,
%
%e4.8 ###
\begin{equation}\label{Eq:setdimUpp1}
\dimh X(E) \le\min\biggl\{d, \frac1 {H_1} \dimh E \biggr\}\qquad
\mathrm{a.s.}
\end{equation}
\end{proposition}

\begin{pf}
Without loss of generality, we assume that $E \subset[0, 1]^N$.
For any constant $\gamma\in(\dimh E, \beta H_1)$, there exists
a sequence of balls $\{B(t_k, h_k), k \ge1\}$ such that
%
%e4.9 ###
\begin{equation}\label{Eq:coverE}
E \subset\limsup_{k \to\infty} B(t_k, h_k) \quad\mbox{and}\quad
\sum_{k=1}^\infty(2 h_k)^\gamma< \infty.
\end{equation}

For a constant $M > 0$, let $\Omega_M = \{\omega\dvt \sup_{t \in[0,
1]^N} \|X(t)\| \le M\}$. Since the sample function of $X(t)$ is
almost surely bounded on $[0, 1]^N$, we have $\lim_{M \to\infty}
\P(\Omega_M) = 1.$ Note that $
X(E) \subset\limsup_{k \to\infty} B(X(t_k), D(t_k, h_k))
$ and, by Lemma \ref{Lem:Momentest}, (\ref{Eq:coverE}) and the fact that
$\gamma< \beta H_1$, we have
\begin{eqnarray}\label{Eq:coverXE2}
\sum_{k=1}^\infty\E(D(t_k, h_k)^{\gamma/H_1}
\I_{\Omega_M} ) &\le& \mathrm{e}^{2M} \sum_{k=1}^\infty\E\bigl(D(t_k,
h_k)^{\gamma/ {H_1}} \mathrm{e}^{-D(t_k, h_k)} \bigr)\nonumber\\[-8pt]\\[-8pt]
&\le& \mathrm{e}^{2M} K_{{4,1}} \sum_{k=1}^\infty h_k^{\gamma} < \infty.\nonumber
\end{eqnarray}
It follows from (\ref{Eq:coverXE2}) that
$\sum_{k=1}^\infty D(t_k, h_k)^{{\gamma}/ {H_1}} < \infty$
almost surely on $\Omega_M.$
This implies that $\dimh X(E) \le
\gamma/H_1$ almost surely on $\Omega_M$. Letting $M \to\infty$
first and then $\gamma\downarrow\dimh E$ along the rational
numbers proves (\ref{Eq:setdimUpp1}).
\end{pf}

Putting Proposition \ref{Prop:LowerBD2} and Proposition
\ref{Prop:setdimUppBD} together, we derive the following theorem.

\begin{theorem}\label{Th:dimXE2}
Let $X=\{X(t), t \in\R^N\}$ be a random field with
values in $\R^d$ whose sample function is a.s. bounded on all compact
subsets of $\R^N$. If there is a constant
$H >0$ such that for every $\ep
> 0$, $X$ satisfies conditions \textup{(C1)} with $H_1 = H-\ep$
and (\textup{C}$2'$) with $H_2 = H + \ep$, then for every analytic set $E
\subset\R^N$ that satisfies $\dimh E < \beta H$,
%
%e4.10 ###
\begin{equation}\label{Eq:dimXE2}
\dimh X(E) = \min\biggl\{d, \frac1 {H} \dimh E \biggr\}\qquad
\mathrm{a.s.}
\end{equation}
\end{theorem}

%In Theorems \ref{Th:main2a} and \ref{Th:dimXE2} the exceptional
%null probability events depend on $E$. It would be interesting
%to investigate whether or when there is a single exceptional
%null probability event $\Omega_0$ such that, for all $\omega
%(\ref{Eq:dimXE2}) hold simultaneously for all Borel
%sets $E \subset\R^N$. Such a ``\emph{uniform} dimension result''
%is applicable even when $E$ is random. For uniform dimension
%results for Brownian motion, stable L\'evy processes and fractional
%Brownian motion, see the survey papers of Pruitt (1975),
%Taylor (1986) and Xiao (2004). The problem is open
%for all the non-Gaussian random fields considered in this paper.

%s5 ###
\section{Applications}
\label{sec:examples}

The general results in Sections \ref{Sec:measures} and \ref{Sec:set} can be applied to wide
classes of Gaussian or non-Gaussian random fields. Since the
applications to Gaussian random fields can be carried out by
extending Xiao (\citeyear{Xiao07}, \citeyear{Xiao09}), we
will focus on non-Gaussian random fields in this section.

%s5.1 ###
\subsection{Self-similar stable random fields}

If $X=\{X(t), t \in\R_+\}$ is a stable L\'evy process in
$\R^d$, the Hausdorff dimensions of its image sets have been
well studied; see Taylor (\citeyear{Taylor86a}) and Xiao (\citeyear{Xiao04}) for historical
accounts. The packing dimension results similar to those
in Sections \ref{Sec:measures} and \ref{Sec:set} have also been obtained
by Khoshnevisan, Schilling and Xiao (\citeyear{KSX09}) for L\'evy processes.
In this subsection, we will only consider non-Markov stable
processes and stable random fields.

Let $X_0=\{X_0(t), t \in\R^N\}$ be an $\alpha$-stable random
field in $\R$ with the representation
%
%e5.1 ###
\begin{equation}\label{Eq:stablerep}
X_0(t) = \int_{F} f(t, x) M(\mathrm{d}x),
\end{equation}
where $M$ is a symmetric $\alpha$-stable (S$\alpha$S) random
measure on a measurable space $(F, \mathcal{F})$ with control
measure $m$ and $f(t, \cdot)\dvtx F \to\R$ ($t \in\R^N$) is a
family of functions on $F$ satisfying
\[
%begin{equation}\label{Eq:fcon1}
\int_{F} |f(t, x)|^\alpha m(\mathrm{d}x) < \infty \qquad\forall t
\in\R^N.
\]
%
%end{equation}
For any integer $n \ge1$ and $t_1, \ldots, t_n \in
\R^N$, the characteristic function of the joint distribution of
$X_0(t_1), \ldots, X_0(t_n)$ is given by
\[
\E\exp\Biggl(\mathrm{i} \sum_{j=1}^n \xi_j X_0(t_j) \Biggr) = \exp\Biggl(- \Biggl\|
\sum_{j=1}^n \xi_j f(t_j, \cdot) \Biggr\|_{\a, m}^\a\Biggr) ,
\]
where $\xi_j \in\R$ ($1 \le j\le n$) and $\|\cdot\|_{\a, m}$ is the
$L^\a(F, \mathcal{F}, m)$-norm (or quasi-norm if $\alpha< 1$).

The class of $\alpha$-stable random fields with representation
(\ref{Eq:stablerep}) is broad. In particular, if a random
field $X_0=\{X_0(t), t \in\R^N\}$ is $\alpha$-stable with
$\alpha\ne1$ or symmetric $\alpha$-stable, and is \emph{separable in
probability} (i.e., there is a countable subset $T_0 \subset
\R^N$ such that for every $t \in\R^N$, there exists a sequence
$\{t_k\} \subset T_0$ such that $X_0(t_k ) \to X_0(t)$ in
probability), then $X_0$ has a representation
(\ref{Eq:stablerep}); see Theorems 13.2.1 and 13.2.2 in
Samorodnitsky and Taqqu (\citeyear{ST94}).

For a separable $\alpha$-stable random field in $\R$ given by
(\ref{Eq:stablerep}), Rosinski and Samorodnitsky (\citeyear{RosinskiSam93})
investigated the asymptotic behavior of $\P\{\sup_{t \in[0,
1]^N} |X_0(t)| \ge u \}$ as $u \to\infty$ (see also
Samorodnitsky and Taqqu (\citeyear{ST94})). The following lemma
is a consequence of their result.
\begin{lemma}\label{Eq:stablesuptail}
Let $X_0= \{X_0(t), t \in\R^N\}$ be a separable $\alpha$-stable
random field in $\R$ given in the form (\ref{Eq:stablerep}).
Assume that $X_0$ has a.s. bounded sample paths on $[0, 1]^N$.
There then exists a positive and finite constant $K_{{5, 1}}$,
depending on $\a$, $f$ and $m$ only, such that for all $u
> 0$,
%
%e5.2 ###
\begin{equation}\label{Eq:RoSam93}
\lim_{u \to\infty} u^{\a} \P\Bigl\{\sup_{t \in[0, 1]^N} |X_0(t)|
\ge u \Bigr\} = K_{{5, 1}} .
\end{equation}
\end{lemma}

\begin{remark}
In the above lemma, it is crucial to assume that $X_0$ has bounded
sample paths on $[0, 1]^N$ almost surely. Otherwise,
(\ref{Eq:RoSam93}) may not hold, as shown by the linear fractional
stable motion $X_0$ with $0 < \a< 1$ (see Example \ref{Ex:lfs}
below).
\end{remark}

We define an $\a$-stable random field $X=\{X(t), t \in\R^N\}$
with values in $\R^d$ by
%
%e5.3 ###
\begin{equation}\label{Def:Xd}
X(t) = (X_1(t), \ldots, X_d(t) ),
\end{equation}
where $X_1, \ldots, X_d$ are independent copies of $X_0$.

The following result gives the Hausdorff and packing dimensions of
the image measures of self-similar stable random fields.
\begin{theorem}\label{Th:Stableimage}
Let $X=\{X(t), t \in\R^N\}$ be a separable $\alpha$-stable
field with values in $\R^d$ defined by (\ref{Def:Xd}), where $X_0$
is given in the form (\ref{Eq:stablerep}). Suppose that $X_0$ is
$H$-SSSI and its sample path is
a.s. bounded on all compact subsets of $\R^N$. Then, for every
finite Borel measure $\mu$ on $\R^N$,
%
%e5.4 ###
\begin{equation}\label{Eq:stablehdim}
\dimh\mu_{X} = \min\biggl\{d, \frac1 {H} \dimh\mu \biggr\}
\quad\mbox{and}\quad
\dimp\mu_{X} = \frac1 {H} \Dim_{Hd} \mu \qquad\mathrm{a.s.}
\end{equation}
Moreover, for every analytic set $E \subset\R^N$ that
satisfies $\dimh E < \alpha H$, we have
\[
\dimh X(E) = \min\biggl\{d, \frac1 {H} \dimh E \biggr\} \qquad\mathrm{a.s.}
\]
\end{theorem}

\begin{pf}
It follows from the self-similarity and Lemma
\ref{Eq:stablesuptail} that $X$ satisfies condition (C1) with $H_1=H$ and
$\beta= \alpha$. On the other hand, condition (C2) with $H_2 = H$ is
satisfied because $X$ is $H$-self-similar and has stationary
increments, and the $\a$-stable variable $X(1)$ has a bounded
continuous density function. Therefore, both equalities in (\ref
{Eq:stablehdim})
follow from Theorems 3.8 and~3.12.
Finally, the last conclusion follows from Theorem \ref{Th:dimXE2}.
\end{pf}

Next, we consider two important types of SSSI
stable processes.

\begin{example}[(Linear fractional stable motion)]\label{Ex:lfs}
Let $0 < \alpha< 2$ and $H \in(0, 1)$ be given constants.
We define an $\a$-stable process $X_0 = \{X_0(t), t \in\R_+\}$
with values in $\R$ by
%
%e5.5 ###
\begin{equation}\label{Eq:Repstable}
X_0(t) = \int_{\R} h_{H} (t, s) M_\a(\mathrm{d}s) ,
\end{equation}
where $M_\a$ is a symmetric $\alpha$-stable random measure on $\R$
with Lebesgue measure as its control measure and where
\[
h_{H} (t, s) = a \{(t-s)_+^{H-1/\a} -
(-s)_+^{H-1/\a} \}
+ b \{(t-s)_{-}^{H-1/\a} - (-s)_{-}^{H-1/\a}
\} .
\]
In the above, $a, b \in\R$ are constants with $|a|+|b| \ne0$,
$t_{+} = \max\{t, 0\}$ and $t_{-} = \max\{-t, 0\}$. The
$\a$-stable process $X_0$ is then $H$-self-similar with stationary
increments, which is called an $(\a, H)$\emph{-linear fractional stable
motion}. If $H = \frac1 {\a}$, then the integral in
(\ref{Eq:Repstable}) is understood as $a M([0, t])$ if $t \ge0$
and as $b M([t, 0])$ if $t < 0$. Hence, $X_0$ is an
$\a$-stable L\'evy process.

Maejima (\citeyear{Maejima83}) proved that if $ \a H < 1$, then $X_0$ is a.s. unbounded
on any interval of positive length. On the other hand, if
$ \a H > 1$ (i.e., $1 <\a< 2$ and $1/\a< H <1$), then
Kolmogorov's continuity theorem implies that $X_0$ is a.s. continuous.
In the latter case, Takashima (\citeyear{Taka89}) further studied
the local and uniform H\"older continuity of $X_0$. His Theorems
3.1 and 3.4 showed that the local H\"older exponent of $X_0$ equals
$H$. However, the exponent of the uniform H\"older continuity cannot be
bigger than $H - \frac1 \a$.

Now, let $X =\{X(t), t \in\R_+\}$ be the $(\a, H)$-linear
fractional stable motion with values in $\R^d$ defined by
(\ref{Def:Xd}). It follows from Theorem \ref{Th:Stableimage}
that if $\a H > 1$, then for every finite Borel measure $\mu$ on $\R_+$,
\[
\dimh\mu_{X} = \min\biggl\{d, \frac1 {H} \dimh\mu\biggr\}
\quad\mbox{and}\quad
\dimp\mu_{X} = \frac1 {H} \Dim_{Hd} \mu \qquad\mathrm{a.s.},
\]
and for every analytic set $E \subset\R_+$, $
\dimh X(E) = \min\{d, \frac1 {H} \dimh E \}$ a.s.
Note that the above dimension results do not depend on the uniform
H\"older exponent of $X$.

There are several ways to define linear fractional $\a$-stable random
fields; see Kokoszka and Taqqu (\citeyear{PKokoszkaTaqqu94}). For example, for $H \in(0,
1)$ and $\a\in(0, 2)$, define
%
%e5.6 ###
\begin{equation}\label{Eq:Linstable}
Z^{H}(t) = \int_{\R^N} (\|t - s\|^{H-  N/ {\alpha}} -
\|s\|^{H-  N /{\alpha}} ) M_\a(\mathrm{d}s)\qquad \forall t
\in\R^N,
\end{equation}
where $M_\a$ is an S$\a$S random measure on $\R^N$ with the
$N$-dimensional Lebesgue measure as its control measure. This is
the stable analog of the $N$-parameter fractional Brownian
motion. However, it follows from Theorem 10.2.3 in Samorodnitsky
and Taqqu (\citeyear{ST94}) that, whenever $N \ge2$, the sample paths of
$Z^{H}$ are a.s. unbounded on any interval in $\R^N$. Thus, the
results of this paper do not apply to $Z^{H}$ when $N \ge2$. In
general, little is known about the sample path properties of
$Z^{H}$.
\end{example}

\begin{example}[(Harmonizable fractional stable motion)]\label{Ex:lhsm}
Given $0 < \alpha< 2$ and $H \in(0, 1)$, the harmonizable fractional
stable field $\widetilde Z^{H} = \{ \widetilde Z^{H}(t), t \in
\R^N\}$ with values in $\R$ is defined by
%
%e5.7 ###
\begin{equation}\label{eq:stableharm}
\widetilde Z^{H} (t) = \Re\int_{\R^N}
\frac{\mathrm{e}^{\mathrm{i} \l t, \lambda\rangle } - 1}{\|\lambda\|^{H+N/\a}}
\widetilde{M}_\a(\mathrm{d}\lambda) ,
\end{equation}
where $\widetilde{M}_\alpha$ is a complex-valued, rotationally invariant
$\a$-stable random
measure on $\R^N$ with the $N$-dimensional Lebesgue measure as its
control measure. It is easy to verify that the $\a$-stable random
field $\widetilde Z^{H}$ is $H$-self-similar with stationary
increments.

It follows from Theorem 10.4.2 in Samorodnitsky and Taqqu (\citeyear{ST94})
(which covers the case $0 < \a< 1$) and Theorem 3 of Nolan (\citeyear{Nolan89})
(which covers $1 \le\alpha< 2$) that $\widetilde Z^{H}$ has\vspace*{1pt}
continuous sample paths almost surely. Moreover, it can be proven
that $\widetilde Z^{H}$ satisfies the following uniform H\"older
continuity: for any compact
interval $I = [a, b] \subset\R^N$ and any $\ep> 0$,
%
%e5.8 ###
\begin{equation}\label{Eq:harmHolder}
\lim_{h \to0} \mathop{\sup_{s, t \in I}}_{
\|s-t\|\le h} \frac{|\widetilde Z^{H}(t) - \widetilde
Z^{H}(s)|} {\|t-s\|^H |\log\|t-s\| |^{1/ 2 + 1/
\a+ \ep}} = 0 \qquad\mathrm{a.s.}
\end{equation}
When $N=1$, (\ref{Eq:harmHolder}) is due to K\^ono and Maejima
(\citeyear{KoMae91}). In general, (\ref{Eq:harmHolder}) follows from the results
in Bierm\'e and Lacaux (\citeyear{BiermeLacaux07}) or Xiao (\citeyear{Xiao09b}). Note that the H\"older
continuity of $\widetilde Z^{H}$ is
different from that of the linear fractional stable motions.

Applying Theorem \ref{Th:main2a} to the harmonizable fractional
stable motion in $\R^d$ defined as in (\ref{Def:Xd}), still denoted
by $\widetilde Z^{H}$, we derive that for every analytic set $E \subset
\R^N$,
%
%e5.9 ###
\begin{equation}\label{Eq:stabledim3}
\dimh\widetilde Z^{H}(E) = \min\biggl\{d, \frac1 {H} \dimh
E \biggr\}\quad \mbox{and}\quad
\dimp\widetilde Z^{H}(E) = \frac1 {H} \Dim_{Hd} E\qquad
\mathrm{a.s.}
\end{equation}
\end{example}

\begin{remark}
The results in this section are applicable
to other self-similar stable random fields, including the Telecom process
(L\'evy and Taqqu (\citeyear{LevyTaqqu00}), Pipiras and Taqqu (\citeyear{PipirasTaqqu00})), self-similar
fields of L\'evy--Chentsov type (Samorodnitsky and Taqqu (\citeyear{ST94}),
Shieh (\citeyear{Sh96})) and the stable sheet (Ehm (\citeyear{Ehm})). We leave the details to
interested readers.
\end{remark}

%s5.2 ###
\subsection{Real harmonizable fractional L\'evy motion}

We show that the results in Sections \ref{Sec:measures} and \ref{Sec:set} can be applied to
the real harmonizable fractional L\'evy motion (RHFLM)
introduced by Benassi, Cohen and Istas (\citeyear{BenassiCohen02}). To recall
their definition, let $\nu$ be a Borel measure on $\C$
which satisfies $\int_{\C} |z|^p \nu(\mathrm{d}z) < \infty$
for all $p \ge2.$ We assume that $\nu$ is rotationally
invariant. Hence, if $P$ is the map $z= \rho \mathrm{e}^{\mathrm{i}\theta}
\mapsto(\theta, \rho) \in[0, 2 \curpi) \times\R_+$,
then the image measure of $\nu$ under $P$ can be written as
$\nu_{P}(\mathrm{d} \theta, \mathrm{d} \rho) = \mathrm{d} \theta\nu_\rho(\mathrm{d}\rho),$
where $\mathrm{d}\theta$ is the uniform measure on $[0, 2 \curpi)$ and
$\nu_\rho$ is a Borel measure on $\R_+$.

Let $N(\mathrm{d} \xi, \mathrm{d}z)$ be a Poisson random measure on $\R^d\times\C$
with mean measure
$n (\mathrm{d}\xi, \mathrm{d}z) = \E(N(\mathrm{d} \xi, \mathrm{d}z))= \mathrm{d} \xi\,\nu(\mathrm{d}z)$ and let
$\widetilde N
(\mathrm{d} \xi, \mathrm{d}z)= N(\mathrm{d} \xi, \mathrm{d}z)
- n(\mathrm{d}\xi, \mathrm{d}z)$ be the compensated Poisson measure. Then, according
to Definition 2.3 in Benassi, Cohen and Istas (\citeyear{BenassiCohen02}),
a real harmonizable fractional L\'evy motion (without the Gaussian part)
$X_0^H = \{X_0^H(t), t \in\R^N\}$ with index $H \in(0, 1)$ is
defined by
%
%e5.10 ###
\begin{equation}\label{Eq:RHFLM}
X_0^H (t) = \int_{\R^N\times\C} 2 \Re\biggl(\frac{\mathrm{e}^{-\mathrm{i} \l t, \xi
\rangle } - 1}
{\|\xi\|^{H+N/2}} z \biggr)
\widetilde{N}(\mathrm{d}\xi, \mathrm{d}z) \qquad\mbox{for all } t \in\R^N.
\end{equation}

As shown by Benassi, Cohen and Istas (\citeyear{BenassiCohen02}),
$X_0^H$ has stationary increments, as well as moments of all orders;
it behaves locally like fractional Brownian motion, but at the large
scale, it behaves like harmonizable fractional stable motion
$\widetilde Z^{H}$ in (\ref{eq:stableharm}).
Because of these multiscale properties, RHFLM's form a class of
flexible stochastic models.

The following equation on characteristic
functions of $X_0^H$ was given by Benassi, Cohen and Istas (\citeyear{BenassiCohen02}):
for all integers $n \ge2$, all $t^1, \ldots, t^n \in\R^N$
and all $u^1, \ldots, u^n\in\R$,
%
%e5.11 ###
\begin{equation}\label{Eq:RHChf}
\E\exp\Biggl( \mathrm{i} \sum_{j=1}^n u^j X_0^H (t^j) \Biggr)
= \exp\biggl(\int_{\R^N\times\C} \bigl[ \mathrm{e}^{f_n(\xi, z)} - 1 -
f_n(\xi, z) \bigr] \,\mathrm{d} \xi\,\nu(\mathrm{d}z) \biggr),
\end{equation}
where
\[
f_n(\xi, z) = \mathrm{i} 2\Re\Biggl(z \sum_{j=1}^n u^j \frac{\mathrm{e}^{-\mathrm{i} \l t^j, \xi
\rangle } - 1}
{\|\xi\|^{H+{N}/{2}}} \Biggr).
\]
In particular, for any $s, t \in\R^N$ and $u\in\R$, (\ref
{Eq:RHChf}) gives
that
%
%e5.12 ###
\begin{equation}\label{Eq:RHChf2}
\E\exp\biggl(\mathrm{i} u \frac{X_0^H (t)- X_0^H (s)} {\|t-s\|^H} \biggr)
= \exp\biggl(- 2 \curpi\int_{\R^N} \psi\biggl(\frac{2 u
(1 - \cos\langle t-s, \xi\rangle)} {\|t-s\|^H
\|\xi\|^{H+{N}/{2}}} \biggr) \,\mathrm{d}\xi \biggr),
\end{equation}
where, for every $x \in\R$, $\psi(x)$ is defined by
$ \psi(x) = \int_0^\infty (1 - \cos(x \rho) ) \nu_\rho(\mathrm{d} \rho)$.
Note that the function $\psi$ is non-negative and continuous.
Moreover, up to a constant, it is the characteristic exponent of
the infinitely divisible law in $\C$ with L\'evy measure $\nu$. %It
%is reasonable to us that the behavior of $\psi$ may determine (to
%a large extent) the properties of RHFLM $X_0^H$. This is the
%motivation of condition (\ref{Con:psi}) below.
For the proof of Theorem \ref{Th:RHFLMimage}, we will make use of
the following fact: there exists a positive constant $K $ such
that
%
%e5.13 ###
\begin{equation}\label{Eq:psi-trun}
\psi(x) \ge K^{-1} x^2 \int_0^{x^{-1}} \rho^2
\nu_\rho(\mathrm{d}\rho) \qquad\mbox{for all } x \in[0, 1].
\end{equation}
This is verified by using the inequality $1 - \cos x \ge K^{-1}
x^2$ for all $x \in[0, 1]$.

\begin{theorem}\label{Th:RHFLMimage}
Let $X^H=\{X^H(t), t \in\R^N\}$ be a separable real harmonizable
fractional L\'evy field in $\R^d$ defined by (\ref{Def:Xd}),
where $X_0^H$ is defined as in (\ref{Eq:RHFLM}). Assume that
$\psi$ satisfies the following condition: there exists a constant
$\delta\in(0, 1]$ such that
%
%e5.14 ###
\begin{equation}\label{Con:psi}
\frac{\psi(ax)} {\psi(x)}\ge a^{\delta} \qquad\mbox{for all }
a \ge1 \mbox{ and } x \in\R.
\end{equation}
Then, for every analytic set $ E \subset\R^N$,
%
%e5.15 ###
\begin{equation}\label{Eq:RHFLMhdim}
\dimh X^H(E) = \min\biggl\{d, \frac1 {H} \dimh E \biggr\}\quad
\mbox{and}\quad
\dimp X^H(E) = \frac1 {H} \Dim_{Hd} E \qquad\mathrm{a.s.}
\end{equation}
\end{theorem}

\begin{pf}
It follows from Proposition 3.3 in Benassi,
Cohen and Istas (\citeyear{BenassiCohen02}) that for every $\ep> 0$, $X^H$
satisfies almost surely a uniform
H\"older condition of order $H - \ep$ on all compact sets of $\R^N$.
Hence, the upper bounds in (\ref{Eq:RHFLMhdim})
follow from Proposition \ref{Prop:setUpp1}.

In order to prove the desired lower bounds in (\ref{Eq:RHFLMhdim}),
by Proposition \ref{Prop:LowerBD2}, it
suffices to show that $X^H$ satisfies condition (C$2'$) with $H_2
= H$. This is done by showing that there exists a positive function
$g \in L^1(\R^d)$ such that for all $s, t \in\R^N$ satisfying
$\|s-t\|\le1$, we have
%
%e5.16 ###
\begin{equation}\label{Eq:g}
\bigl|\E\bigl( \mathrm{e}^{\mathrm{i} \l u, X(t) - X(s)\rangle
/\|t-s\|^{H}} \bigr) \bigr|\le g(u) \quad\mbox{for all } u \in
\R^d .
\end{equation}
This and the Fourier inversion formula together imply that the density functions
of $X(t) - X(s)/\|t-s\|^{H}$ are uniformly bounded for all $s, t \in
\R^N$ satisfying $\|s-t\|\le1$.

Since the coordinate processes $X_1^H, \ldots,
X_d^H$ are independent copies of $X_0^H$, it is sufficient to
prove (\ref{Eq:g}) for $d = 1$.
Note that, by (\ref{Eq:g}), we can take $g(u) = 1$ for all $|u | \le
1$. For
any $u$ such that $|u| > 1$, condition (\ref{Con:psi}) implies
that
\begin{eqnarray}\label{Eq:psi-int}
\hspace*{-28pt}\int_{\R^N} \psi\biggl(\frac{2 u (1 - \cos\langle t-s, \xi\rangle)}
{\|t-s\|^H \|\xi\|^{H+N/2}} \biggr) \,\mathrm{d}\xi&\ge& K
|u|^\delta\int_{\R^N} \psi\biggl(\frac{1 - \cos\langle t-s,
\xi\rangle} {\|t-s\|^H \|\xi\|^{H+N/2}} \biggr)
\,\mathrm{d}\xi\nonumber\\[-8pt]\\[-8pt]
&\ge& K |u|^\delta\int_{\|\xi\| \ge\gamma\|t-s\|^{-1}}
\psi\biggl(\frac{1 - \cos\langle t-s, \xi\rangle} {\|t-s\|^H
\|\xi\|^{H+{N}/{2}}} \biggr) \,\mathrm{d}\xi,\nonumber
\end{eqnarray}
where $\gamma> 1$ is a constant whose value will be chosen
later.

By a change of variable $\xi\mapsto\eta \|t-s\|^{-1}$, we see
that the last integral becomes
\begin{eqnarray}\label{Eq:psi-int2}
&& \int_{\|\eta\| \ge\gamma}
\psi\biggl(\frac{\|t-s\|^{N/2}(1 - \cos\langle
(t-s)/\|t-s\|, \eta\rangle)} { \|\eta\|^{H+N/2}} \biggr)
\frac{\mathrm{d}\eta} {\|t-s\|^N}\nonumber\\[-8pt]\\[-8pt]
&&\quad \ge K \int_{\|\eta\| \ge\gamma} \frac{(1 - \cos\langle
(t-s)/\|t-s\|, \eta\rangle)^2} { \|\eta\|^{2H+N}} \,\mathrm{d}\eta,\nonumber
\end{eqnarray}
where the inequality follows from (\ref{Eq:psi-trun}), and we have
used the fact that $\|t-s\| \le1$ and taken $\gamma$ large. The
last integral is a constant because the Lebesgue measure is
rotationally invariant. Thus, we have proven that for $|u| > 1$,
%
%e5.17 ###
\begin{equation}\label{Eq:RHChf3}
\E\exp\biggl(\mathrm{i} u \frac{X_0^H (t)- X_0^H (s)} {\|t-s\|^H} \biggr) \le
\exp(-K_{{5,2}} |u|^\delta) .
\end{equation}
Therefore, when $d = 1$, (\ref{Eq:g}) holds for the function $g$
defined as
$g(u) = 1$ if $|u| \le1$ and $g(u)=
\mathrm{e}^{-K_{{5,2}} |u|^\delta}$ if $|u| > 1.$
This completes the proof of Theorem \ref{Th:RHFLMimage}.
\end{pf}

We mention that Benassi, Cohen and Istas (\citeyear{BenassiCohen04}) have introduced another
interesting class of fractional L\'evy fields, namely, the moving average
fractional L\'evy fields (MAFLF). Similarly to the contrast between
linear fractional stable motion and harmonizable fractional stable motion,
many properties of MAFLF's are different from those of RHFLM's.
For example, the exponent of the uniform modulus
of continuity of an MAFLF is strictly smaller than its local H\"older
exponent. Nevertheless, we believe that the arguments
in this paper are applicable to MAFLF's. This and some related problems
will be dealt with elsewhere.

%s5.3 ###
\subsection{The Rosenblatt process}

Given an integer $m \ge2$ and a constant $\kappa\in(1/2-
1/(2m), 1/2)$, the Hermite process $Y^{m, \kappa}= \{Y^{m,
\kappa}(t), t \in\R_+\}$ of order $m$ is defined by
%
%e5.18 ###
\begin{equation}\label{Eq:Hermite}
Y^{m, \kappa}(t) = K_{{5,3}} \int_{\R^m}' \Biggl\{\int_0^t
\prod_{j=1}^m (s - u_j)_+^{\kappa- 1} \,\mathrm{d}s \Biggr\}
\,\mathrm{d}B(u_1)\,\cdots\,
\mathrm{d}B(u_m) ,
\end{equation}
where $K_{{5,3}} > 0$ is a normalizing constant depending on $m$
and $\kappa$ only and the integral $\int_{\R^m}'$ is the $m$-tuple
Wiener--It\^o integral with respect to the standard Brownian
motion excluding the diagonals $\{u_i = u_j\}$, $i \ne j$.
The integral (\ref{Eq:Hermite}) is also well defined if $m=1$;
the process is a fractional Brownian motion for which
the problem considered in this paper has been solved.

The Hermite process $Y^{m, \kappa}$ is $H$-SSSI and $H
= 1 + m \kappa- \frac m 2 \in
(0, 1)$. It is a non-Gaussian process and often appears in non-central
limit theorems for processes defined as integrals or partial
sums of nonlinear functionals of stationary Gaussian
sequences with long-range dependence; see Taqqu (\citeyear{Taqqu75}, \citeyear{Taqqu79}),
Dobrushin and Major (\citeyear{DMajor79}) and Major (\citeyear{Major}).

It follows from Theorem 6.3 of Taqqu (\citeyear{Taqqu79}) that the Hermite
process $Y^{m, \kappa}$ has the following equivalent
representation:
%
%e5.19 ###
\begin{equation}\label{Eq:HermiteHarm}
Y^{m, \kappa}(t)= K_{{5,4}} \int^{\prime}_{\R^{m}}
\frac{\mathrm{e}^{\mathrm{i}t(u_1 +\cdots+u_m)} -1} {\mathrm{i}(u_1 + \cdots+u_m)}
\prod_{j=1}^m |u_{j}|^{\kappa-1} Z_G(\mathrm{d}u_{1})\cdots Z_G(\mathrm{d}u_{m}) ,
\end{equation}
where $K_{{5,4}} > 0$ is a normalizing constant and $Z_G$ is a
centered complex Gaussian random measure on $\R$ with Lebesgue
measure as its control measure.

%Let $0<\kappa<N/2m$, and let $\phi_t(u)$ be
%$$\phi_t(u)=\prod_{j=1}^N \frac{e^{t_ju_j} -1}{iu_j},\ \
%t=(t_1,\cdots,t_N),u=(u_1,\cdots,u_N).$$
%The $N-$parameter {\it Hermite process of rank m$\geq2$} is
%defined to be
%$$
%Y(t)=\int^{\prime\prime}_{R^{Nm}}
%|u^{(m)}|^{\kappa-N/2}Z_G(du^{(1)})\cdots Z_G(du^{(m)}),
%$$
%where the integral in the above means $m-$fold Wiener--It\^o integral.
% note that the subs and ups notation is reversed here and in Major's

Mori and Oodaira (\citeyear{MO86}) studied the functional laws of the
iterated logarithms for the Hermite process $Y^{m, \kappa}$. Lemma
\ref{MWtail} follows from Lemma 6.3 in Mori and Oodaira
(\citeyear{MO86}).
\begin{lemma}\label{MWtail}
There exist positive constants $K_{{5, 5}}$ and $K_{{5, 6}}$,
depending on $m$ only, such that
$
\P\{ \max_{t \in[0, 1]} |Y^{m, \kappa}(t)| \ge u \} \le
\exp(-K_{{5, 6}} u^{2/m} )$ for all $u \ge K_{{5, 5}}.
$
\end{lemma}

Using Lemma \ref{MWtail}, one can derive easily a uniform
modulus of continuity for $Y^{m, \kappa}$.

\begin{lemma}\label{Lem:modcon}
There exists a finite constant $K_{{5, 7}}$
such that for all constants $0\le a<b <\infty$,
%
%e5.20 ###
\begin{equation}\label{Eq:uniform}
\limsup_{h\downarrow0} \sup_{a \leq t\leq b -h}\sup_{0\leq
s\leq h} \frac{ |Y^{m, \kappa}(t+s)-Y^{m, \kappa}(t) |}{h^{H} (
\log1/h )^{m/2}}\leq K_{{5, 7}}\qquad
\mathrm{a.s.},
\end{equation}
where $H = 1 + m \kappa- \frac m 2$.
\end{lemma}

\begin{pf}
For every $t \ge0$ and $h > 0$, the self-similarity of
$Y^{m, \kappa}$ and Lemma \ref{MWtail} together imply that
%
%e5.21 ###
\begin{equation}\label{Eq:Lamom4}
\P\{|Y^{m, \kappa}(t+h)-Y^{m, \kappa}(h)| > h^{H} u \} \leq
\exp (- K_{{5, 6}} u^{2/m} ).
\end{equation}
Hence, $Y^{m, \kappa} = \{Y^{m, \kappa}(t), t \ge0\}$ satisfies the
conditions of
Lemmas 2.1 and 2.2 in Cs\'aki and Cs\"org\H o (\citeyear{CaskiC}) with $\sigma(h) =
h^{H}$ and $\beta= 2/m$. Consequently, (\ref{Eq:uniform}) follows
directly from Theorem 3.1 in Cs\'aki and Cs\"org\H o (\citeyear{CaskiC}).
\end{pf}

The case $m=2$ has recently received considerable attention. The
process $Y^{2, \kappa}$ is called the \emph{Rosenblatt process} by Taqqu
(\citeyear{Taqqu75}) (or \textit{fractional Rosenblatt motion} by Pipiras
(\citeyear{Pipiras04})). Its self-similarity index is given by $H = 2
\kappa$. This non-Gaussian process in many ways resembles
fractional Brownian motion. For example, since $H > 1/2$,
fractional noise of $Y^{2, \kappa}$ exhibits long-range
dependence. Besides its connections to non-central limit theorems, the
Rosenblatt process also appears in limit theorems for some quadratic
forms of random variables with long-range dependence. Albin (\citeyear{Albin98a},
\citeyear{Albin98b}) has discussed distributional properties and the extreme value
theory of $Y^{2, \kappa}$. In particular, Albin (\citeyear{Albin98b}), Section 16,
obtained sharp asymptotics on the tail probability of $\max_{t \in[0,
1]}Y^{2, \kappa}(t)$. Pipiras (\citeyear{Pipiras04}) established a wavelet-type
expansion for the
Rosenblatt process. Tudor (\citeyear{Tudor}) has recently developed a
stochastic calculus for $Y^{2, \kappa}$ based on both pathwise type
calculus and Malliavin calculus.

We now consider the Rosenblatt process $X^{2, \kappa}$ with values
in $\R^d$ by letting its component processes be independent copies
of $Y^{2, \kappa}$. The following result determines the
Hausdorff and packing dimensions of the image sets of $X^{2, \kappa}$.

\begin{corollary} %\label{Th:Hermite}
Let $X^{2, \kappa} = \{X^{2, \kappa}(t), t \in\R_+\}$ be a
Rosenblatt process in $\R^d$ as defined above. Then, for every
analytic set $E \subset\R_+$, we have
%
%e5.22 ###
\begin{equation} \label{Eq:setFRm1}
\dimh X^{2, \kappa}(E) = \min\biggl\{d, \frac1 {2 \kappa} \dimh
E \biggr\} \quad\mbox{and}\quad
\dimp X^{2, \kappa}(E) = \frac1 {2 \kappa} \Dim_{{2 \kappa
d}} E \qquad\mathrm{a.s.}
\end{equation}
\end{corollary}

\begin{pf}
By Lemma \ref{Lem:modcon}, for any $\ep> 0$, $X^{2, \kappa}$
satisfies a uniform
H\"older condition of order $H-\ep$ (where $H= 2 \kappa$)
on all compact intervals in $\R_+$. On the
other hand, it is known that the random variable $Y^{2,
\kappa}(1)$ has a bounded and continuous density (see Davydov (\citeyear{Davydov90}) or
Albin (\citeyear{Albin98a})). Thus, $X^{2, \kappa}$ also satisfies condition (C2)
with $H_2 = 2 \kappa$. Therefore, the two equalities in (\ref{Eq:setFRm1})
follow from Theorem \ref{Th:main2a}.
\end{pf}

\section*{Acknowledgements}
The authors thank Professor Davar Khoshnevisan for stimulating
discussions on packing dimension profiles, Professors Patrik Albin,
David Nualart and Frederi Viens for very helpful suggestions on the
Rosenblatt process and the anonymous referee for comments and
suggestions which have led to significant improvement of the
manuscript. Research of N.-R. Shieh was partially supported by Taiwan
NSC Grant 962115M002005MY3. Research of Y. Xiao was partially
supported by the NSF Grant DMS-0706728 and the National Natural Science
Foundation of China (No. 70871050).

\printhistory

\end{document}